\documentclass[11pt]{article}
\usepackage{}

\usepackage[total={6in, 8.5in}]{geometry}
\usepackage{amsfonts}
\usepackage{amsthm}
\usepackage{amssymb}
\usepackage{amsmath}
\usepackage{enumerate}
\usepackage[
pdfauthor={ESYZ},
pdftitle={Toughness and spanning trees in K4mf graphs},
pdfstartview=XYZ,
bookmarks=true,
colorlinks=true,
linkcolor=blue,
urlcolor=blue,
citecolor=blue,
bookmarks=false,
linktocpage=true,
hyperindex=true
]{hyperref}

\makeatletter
\newcommand{\setword}[2]{%
	\phantomsection
	#1\def\@currentlabel{\unexpanded{#1}}\label{#2}%
}
\makeatother

\makeatletter
\renewcommand*\env@matrix[1][*\c@MaxMatrixCols c]{%
	\hskip -\arraycolsep
	\let\@ifnextchar\new@ifnextchar
	\array{#1}}
\makeatother

\usepackage{graphicx}
\usepackage[natural]{xcolor}
\usepackage{tikz}
\usepackage{tikz-3dplot}
\usetikzlibrary{shapes}
\usetikzlibrary{arrows,decorations.pathmorphing,backgrounds,positioning,fit,petri,automata}
\usetikzlibrary {positioning}
\usepackage{tkz-graph}
\usetikzlibrary{arrows}
\usepackage{tkz-euclide}
\usetkzobj{all}

\usepackage{pgf,tikz,pgfplots}

\usepackage{mathrsfs}

\usetikzlibrary{arrows}

\usepackage{xcolor}
\long\def\ignore#1{}

\let\oldi\ignore

 \oldi{
\newtheorem{THM}{\textbf{Theorem}}[section]
\newtheorem{THMs}{\textbf{Theorem}}[section]
\newtheorem{DEF}[THM]{\textbf{Definition}}[section]
\newtheorem{LEM}[THM]{\textbf{Lemma}}
\newtheorem{CON}[THM]{\textbf{Conjecture}}
\newtheorem{PROP}[THM]{\textbf{Proposition}}
\newtheorem{COR}[THM]{\textbf{Corollary}}
\newtheorem{CORs}{\textbf{Corollary}}[section]
\newtheorem{PRO}[THM]{\textbf{Problem}}
\newcommand{\pf}{\textbf{Proof}.\quad}

\newtheorem{FAC}{\textbf{Fact}}
\newtheorem{REM}{\textbf{Remark}}
\newtheorem{OPR}{\textbf{Operation}}
\newtheorem{CLA}{\textbf{Claim}}[section]
\setcounter{CLA}{1}
 }

\newtheorem{THM}{Theorem}[section]

\newtheorem{LEM}[THM]{Lemma}

\newtheorem{COR}[THM]{Corollary}

\newtheorem{CLA}{Claim}[section]
\newcommand{\pf}{\textbf{Proof}.\quad}

\newtheorem*{LEM1}{\textbf{Lemma 3.1}}
\newtheorem*{LEM2}{\textbf{Lemma 3.2}}
\newtheorem*{LEM3}{\textbf{Lemma 3.3}}
\newtheorem*{LEM4}{\textbf{Lemma 3.4}}

\usepackage{thmtools}
\usepackage{thm-restate}

\newcommand{\CC}{\mathcal{C}}

\newcommand{\pbar}{\overline{\varphi}}
\DeclareMathOperator{\dist}{dist} 

\begin{document}
\title{$\Delta$-critical graphs with a vertex of degree 2}

\author{%
 Yan Cao
 \qquad Guantao Chen\thanks{This author was supported in part by NSF grant DMS-1855716.}\\
  Department of Mathematics and Statistics, \\
  Georgia State University, Atlanta, GA 30302, USA\\
   \texttt{ycao17@gsu.edu}
  \qquad
   \texttt{gchen@gsu.edu}%
 \and 
 Songling Shan \\
 Department of Mathematics, \\
 Illinois State Univeristy, Normal, IL 61790, USA \\
 \texttt{sshan12@ilstu.edu}
} 

\date{May 25, 2020}
\maketitle

 \begin{abstract}
 Let $G$ be a simple graph with maximum degree $\Delta$.
 A classic result of Vizing shows that $\chi'(G)$, the chromatic index of $G$, is either $\Delta$
 or $\Delta+1$.   We say $G$ is of \emph{Class 1} if $\chi'(G)=\Delta$, and is  of \emph{Class 2}
 otherwise.     
 A  graph $G$  is \emph{$\Delta$-critical} if  $\chi'(G)=\Delta+1$ and $\chi'(H)<\Delta+1$ for every proper subgraph $H$ of $G$,
 and is \emph{overfull}  if $|E(G)|>\Delta \lfloor (|V(G)|-1)/2 \rfloor$. 
 Clearly,  overfull graphs are  Class 2. 
 Hilton and Zhao in  1997 conjectured that  if $G$ is obtained 
 from an $n$-vertex $\Delta$-regular Class 1 graph  with  maximum degree greater than $n/3$ by splitting a vertex, 
 then being overfull is the only reason for $G$ to be Class 2. 
 This conjecture was only confirmed  when $\Delta\ge \frac{n}{2}(\sqrt{7}-1)\approx 0.82n$. 
 In this paper, we improve the bound on $\Delta$ from $\frac{n}{2}(\sqrt{7}-1)$ to $0.75n$. 
 Considering the structure of $\Delta$-critical graphs with a vertex of degree 2, 
 we also show that for  an $n$-vertex $\Delta$-critical graph  with  $\Delta\ge \frac{3n}{4}$, 
if it contains a vertex of degree 2, then it is overfull. We actually obtain a more general 
 form of this result, which partially supports the overfull 
 conjecture of  
 Chetwynd and  Hilton from  1986, which states that   if $G$
 is an $n$-vertex $\Delta$-critical graph with $\Delta>n/3$, 
 then $G$ contains an overfull subgraph $H$ with $\Delta(H)=\Delta$. 
Our proof techniques are new and might shed some light on attacking both of the  conjectures 
when $\Delta$ is large. 


 \end{abstract}


\section{Introduction}

For two integer $p,q$ with $q\ge p$, we use $[p,q]$ to denote the set of all integer between $p$ 
and $q$, inclusively. 
We consider only  simple graphs.   Let $G$ be a graph with maximum degree $\Delta(G)=\Delta$.  We denote by $V(G)$ and $E(G)$ the vertex set and edge set of $G$, respectively. 
Let $v\in V(G)$ with $d_G(v)=t\ge 2$ and $N_G(v)=\{u_1,\ldots, u_t\}$.   
A \emph{vertex-splitting} in $G$ at $v$ into two vertices $v_1$ and $v_2$  gives a new graph $G'$ such that 
$V(G')=(V(G)\setminus\{v\})\cup \{v_1,v_2\}$ and $E(G')=E(G-v)\cup \{v_1v_2\}\cup \{v_1u_i: i\in [1,s]\} \cup \{v_2u_i: i\in [s+1,t]\}$, where $s\in [1,t-1]$ is any integer.  We say $G'$ is obtained from $G$ by a vertex-splitting.

An {\it edge $k$-coloring\/} or simply a \emph{$k$-coloring} of $G$ is a mapping $\varphi$ from $E(G)$ to the set of integers
$[1,k]$, called {\it colors\/}, such that  no two adjacent edges receive the same color with respect to $\varphi$.  
The {\it chromatic index\/} of $G$, denoted $\chi'(G)$, is defined to be the smallest integer $k$ so that $G$ has an edge $k$-coloring.  
We denote by $\CC^k(G)$ the set of all edge $k$-colorings of $G$.  
A  graph $G$  is \emph{$\Delta$-critical} if  $\chi'(G)=\Delta+1$ and $\chi'(H)<\Delta+1$ for every proper subgraph $H$ of $G$. 
In 1960's, Vizing~\cite{Vizing-2-classes} showed that a graph of maximum degree
$\Delta$ has  chromatic index either $\Delta$ or $\Delta+1$.
If $\chi'(G)=\Delta$, then $G$ is said to be of {\it Class 1\/}; otherwise, it is said to be
of {\it Class 2\/}.  
Holyer~\cite{Holyer} showed that it is NP-complete to determine whether an arbitrary graph is of Class 1.  
Nevertheless, if a graph $G$ has too many edges, i.e., $|E(G)|>\Delta \lfloor |V(G)|/2\rfloor$,  then we have to color $E(G)$ using exactly  $(\Delta+1)$ colors. Such graphs are called \emph{overfull}. 
Clearly, overfull graphs have an odd order, all regular graphs of an odd order are overfull,
and all graphs obtained from a regular Class 1 graph by a vertex-splitting is overfull.   

Being overfull is definitely a cause for a graph to be Class 2, but is that the only cause? 
Hilton and Zhao~\cite{MR1460574} in 1997 conjectured the following: 
	Let $G$ be an $n$-vertex Class 1 $\Delta$-regular graph with $\Delta>\frac{n}{3}$. If $G^*$
	is obtained from $G$ by a vertex-splitting, then $G^*$ is $\Delta$-critical (vertex-splitting conjecture). 
Clearly, $G^*$ is overfull. The conjecture asserts that for every $e\in E(G^*)$, $G^*-e$
is no longer Class 2. In other words, being overfull is the only cause for $G^*$ to be Class 2. 
This conjecture was verified when $\Delta\ge \frac{n}{2}(\sqrt{7}-1)\approx 0.82n$ by 
Hilton and Zhao~\cite{MR1460574} in 1997.  Song~\cite{MR1874750} in 2002 showed that the conjecture holds for 
a particular class of $n$-vertex Class 1 $\Delta$-regular  graphs with $\Delta\ge \frac{n}{2}$. 
No other progress on this conjecture has been achieved since then. 
We support this conjecture as below.  

\begin{THM}\label{thm:vertex spliting}
	Let $n$ and $ \Delta$ be positive integers such that $\Delta\ge \frac{3(n-1)}{4}$. 
	If $G$ is obtained from an $(n-1)$-vertex $\Delta$-regular Class 1 graph by a vertex-splitting, then $G$ 
	is $\Delta$-critical. 
\end{THM}

If $G$ is an $n$-vertex overfull graph, then $|E(G)|\ge \Delta(n-1)/2 +1$. 
Thus $\sum_{v\in V(G)} (\Delta-d_G(v)) \le \Delta-2$. Therefor, if $G$ has a vertex of degree 2,
then all other vertices of $G$ are of maximum degree. Is the converse of this statement true?
That is, when will be a Class 2 graph with a degree 2 vertex overfull? We investigate this question and show that this happens 
when $\Delta$ is large. In general, for two adjacent vertices $u,v\in V(G)$, 
we call $(u,v)$ a \emph{full-deficiency pair} of $G$
if $d(u)+d(v)=\Delta(G)+2$. In particular, if $v$ is of degree 2 in a $\Delta$-critical graph $G$, then each neighbor  $u$ of 
 $v$ has degree $\Delta$ by Vizing's Adjacency Lemma (this lemma will be introduced in Section 2).  Therefore, $(u,v)$ is a  full-deficiency pair of $G$. 
We obtain the following result. 

\begin{THM}\label{thm:Delta-critical}
	Let $n$ and $ \Delta$ be positive integers such that $\Delta\ge \frac{3(n-1)}{4}$, and  
	$G$ be an  $n$-vertex $\Delta$-critical graph.
	If  $G$ has a full-deficiency pair, then $G$ is overfull. 
	Consequently, $G$ is obtained from an $(n-1)$-vertex $\Delta$-regular Class 1 multigraph by a vertex-splitting.
\end{THM}

Theorem~\ref{thm:Delta-critical} partially supports a conjecture of  
 Chetwynd and  Hilton from 1986~\cite{MR848854,MR975994}. The conjecture states 
the following: Let $G$ be a simple graph  with $\Delta(G)>\frac{1}{3}|V(G)|$. Then $G$ is Class 2 implies that 
$G$ contains an overfull subgraph $H$  with $\Delta(H)=\Delta(G)$ (overfull conjecture). 
The overfull conjecture was confirmed only for some special classes of graphs. 
Chetwynd and   Hilton~\cite{MR975994} in 1989 verified the conjecture for 
$n$-vertex graphs with $\Delta\ge n-3$.
 Hoffman and  Rodger~\cite{comm} in 1992 confirmed the conjecture for complete multipartite graphs. 
Plantholt~\cite{MR2082738} in 2004 showed that the overfull conjecture is affirmative  for 
graphs $G$ with an even order $n$, maximum degree $\Delta$ and minimum degree $\delta$ satisfying 
$(3\delta-\Delta)/2 \ge cn $ for any $c\ge \frac{3}{4}$. 
The overfull conjecture was also confirmed for large regular graphs in 2013~\cite{MR3185848,MR3545109}.

Both the overfull conjecture and the vertex-splitting conjecture are best possible in terms of the 
condition on the maximum degree,  by considering the critical Class 2  graph $P^*$, which is obtained from the Petersen graph by deleting a vertex.   Hilton and Zhao~\cite{MR1460574} proved that the overfull conjecture implies the 
vertex-splitting conjecture. 
 
The results in Theorem~\ref{thm:vertex spliting} and Theorem~\ref{thm:Delta-critical} together 
imply that all $n$-vertex $\Delta$-critical graphs 
with a vertex of degree 2 can be obtained from an $(n-1)$-vertex $\Delta$-regular Class 1 multigraph 
by a vertex splitting, when $\Delta\ge \frac{3(n-1)}{4}$. 
Thereby, these results provide  
 a way of constructing dense $\Delta$-critical graphs, 
which 
are known to be hard. 
The reminder of this paper is organized as follows. 
We introduce some definitions and preliminary results in Section 2. 
In Section 3, we prove Theorem~\ref{thm:vertex spliting} and Theorem~\ref{thm:Delta-critical} 
 by assuming the truth of 
several lemmas. These lemmas will be proved   in the last section.

\section{Definitions and Preliminary Results}\label{lemma}

Let $G$ be a graph. 
For $e\in E(G)$, $G-e$
denotes the graph obtained from $G$ by deleting the edge $e$. 
The symbol $\Delta$  is reserved for $\Delta(G)$, the maximum degree of $G$
throughout  this paper.  A \emph{$k$-vertex}  in $G$  is a vertex of degree exactly $k$
in $G$, and a \emph{$k$-neighbor}  of a vertex $v$ is a neighbor of $v$ that is a $k$-vertex in $G$.
For $u,v\in V(G)$, we use $\dist_G(u,v)$ to denote the distance between $u$ and $v$, which is the length of a shortest path connecting $u$
and $v$ in $G$. For $S\subseteq V(G)$, define $\dist_G(u,S)=\min_{v\in S} \dist_G(u,v)$. 

An edge $e\in E(G)$ is a \emph{critical edge} of $G$ if $\chi'(G-e)<\chi'(G)$. 
It is not hard to see that if 
 $G$ is connected,  $\chi'(G)=\Delta+1$ and every edge of $G$ is critical, then $G$ is $\Delta$-critical.  
Critical graphs are useful since they provide more information about the  structure around a vertex than general Class 2 graphs. For 
example,  Vizing's Adjacency Lemma (VAL) from 1965~\cite{Vizing-2-classes} is a useful tool that reveals certain structure at a vertex
by assuming the criticality of an edge.

\begin{LEM}[Vizing's Adjacency Lemma (VAL)]Let $G$ be a Class 2 graph with maximum degree $\Delta$. If $e=xy$ is a critical edge of $G$, then $x$ has at least $\Delta-d_G(y)+1$ $\Delta$-neighbors in $V(G)\setminus \{y\}$.
	\label{thm:val}
\end{LEM}
Let $G$ be a graph and 
$\varphi\in \CC^k(G-e)$ for some edge $e\in E(G)$ and some integer $k\ge 0$. 
For any $v\in V(G)$, the set of colors \emph{present} at $v$ is 
$\varphi(v)=\{\varphi(f)\,:\, \text{$f$ is incident to $v$}\}$, and the set of colors \emph{missing} at $v$ is $\pbar(v)=[1,k]\setminus\varphi(v)$.  
For a vertex set $X\subseteq V(G)$,  define 
$$
\pbar(X)=\bigcup _{v\in X} \pbar(v).
$$
The set $X$ is called \emph{elementary} with respect to $\varphi$  or simply \emph{$\varphi$-elementary} if $\pbar(u)\cap \pbar(v)=\emptyset$
for every two distinct vertices $u,v\in X$.   Sometimes, we just say that $X$ 
is elementary if the  edge coloring is understood.  





  For two distinct colors $\alpha,\beta \in [1,k]$, let  $H$ be the subgraph of $G$
with $V(H)=V(G)$ and $E(H)$ consisting of edges from $E(G)$ that are colored by $\alpha$
or $\beta$ with respect to $\varphi$. Each component of $H$ is either 
an even cycle or a path, which is called an \emph{$(\alpha,\beta)$-chain} of $G$
with respect to $\varphi$.  If we interchange the colors $\alpha$ and $\beta$
on an $(\alpha,\beta)$-chain $C$ of $G$, we get a new edge $k$-coloring  of $G$, 
and we write $$\varphi'=\varphi/C.$$
This operation is called a \emph{Kempe change}. 
For a color $\alpha$, a sequence of 
{\it Kempe  $(\alpha,*)$-changes}  is a sequence of  
Kempe changes that each involves the exchanging of the color $\alpha$
and another color from $[1,k]$.

Let  $x,y\in V(G)$, and  $\alpha, \beta, \gamma\in [1,k]$ be three colors.   If $x$ and $y$
are contained in a same  $(\alpha,\beta)$-chain of $G$ with respect to $\varphi$, we say $x$ 
and $y$ are \emph{$(\alpha,\beta)$-linked} with respect to $\varphi$.
Otherwise, $x$ and $y$ are \emph{$(\alpha,\beta)$-unlinked} with respect to $\varphi$. Without specifying $\varphi$, when we just say  $x$ and $y$ are $(\alpha,\beta)$-linked or $x$ and $y$ are $(\alpha,\beta)$-unlinked, we mean they are linked or unlinked with respect to the current edge coloring. 
Let $P$ be an 
$(\alpha,\beta)$-chain of $G$ with respect to $\varphi$ that contains both $x$ and $y$. 
If $P$ is a path, denote by $\mathit{P_{[x,y]}(\alpha,\beta, \varphi)}$  the subchain  of $P$ that has endvertices $x$
and $y$.  By \emph{swapping  colors} along  $P_{[x,y]}(\alpha,\beta,\varphi)$, we mean 
exchanging the two colors $\alpha$
and $\beta$ on the path $P_{[x,y]}(\alpha,\beta,\varphi)$. 
The notion $P_{[x,y]}(\alpha,\beta)$ always represents the $(\alpha,\beta)$-chain
with respect to the current edge coloring. 
Define  $P_x(\alpha,\beta,\varphi)$ to be an $(\alpha,\beta)$-chain or an $(\alpha,\beta)$-subchain of $G$ with respect to $\varphi$ that starts at $x$  and ends at a different vertex missing exactly one of $\alpha$ and $\beta$.    
(If $x$ is an endvertex of the $(\alpha,\beta)$-chain that contains $x$, then $P_x(\alpha,\beta,\varphi)$ is unique.  Otherwise, we take one segment of the whole chain to be 
$P_x(\alpha,\beta,\varphi)$. We will specify the segment when it is used.) 
If  $u$  is a vertrex on  $P_x(\alpha,\beta,\varphi)$, we  write  {$\mathit {u\in P_x(\alpha,\beta, \varphi)}$}; and if $uv$  is an edge on  $P_x(\alpha,\beta,\varphi)$, we  write  {$\mathit {uv\in P_x(\alpha,\beta, \varphi)}$}.  Similarly, the notion $P_x(\alpha,\beta)$ always represents the $(\alpha,\beta)$-chain
with respect to the current edge coloring. 
If $u,v\in P_x(\alpha,\beta)$ such that $u$ lies between $x$ and $v$, 
then we say that $P_x(\alpha,\beta)$ \emph{meets $u$ before $v$}. 
Suppose that   $\alpha\in \pbar(x)$ and  $\beta,\gamma\in \varphi(x)$. An $\mathit{(\alpha,\beta)-(\beta,\gamma)}$
\emph{swap at $x$}  consists of two operations:  first swaps colors on $P_x(\alpha,\beta, \varphi)$ to get an edge  $k$-coloring $\varphi'$, and then swaps
colors on $P_x(\beta,\gamma, \varphi')$. 
By convention, an	$(\alpha,\alpha)$-swap at $x$ does  nothing at $x$. 
Suppose the current color of an  edge $uv$ of $G$
is $\alpha$, the notation  $\mathit{uv: \alpha\rightarrow \beta}$  means to recolor  the edge  $uv$ using the color $\beta$. Recall that $\pbar(x)$ is the set of colors not present 
at $x$. 
If $|\pbar(x)|=1$, we will also use $\pbar(x)$ to denote the  color that is missing at $x$. 



Let $\alpha, \beta, \gamma, \tau,\eta\in[1,k]$. 
We will use a  matrix with two rows to denote a sequence of operations  taken  on $\varphi$.
Each entry in the first row represents a path or  a sequence of vertices. 
Each entry in the second row, indicates the action taken on the object above this entry. 
We require the operations to be taken to follow the ``left to right'' order as they appear in 
the matrix. 
For example,   the matrix below indicates 
three sequential operations taken on the graph based 
on the coloring from the previous step:
\[
\begin{bmatrix}
P_{[a, b]}(\alpha, \beta) &   rs & ab \\
\alpha/\beta & \gamma \rightarrow \tau & \eta
\end{bmatrix}.
\]
\begin{enumerate}[Step 1]
	\item Swap colors on the $(\alpha,\beta)$-subchain $P_{[a, b]}(\alpha, \beta,\varphi) $.
	
	\item Do  $rs: \gamma \rightarrow \tau $. 
	\item Color the edge $ab$ using color $\eta$. 
\end{enumerate}

Let   
$T$ be  a sequence of vertices  and edges of  $G$. We denote by \emph{$V(T)$}  
the set of vertices from $V(G)$ that are contained in $T$, and by  
\emph{$E(T)$}  the set of edges  from $E(G)$ that are contained in $T$.

%
%

\subsection{Multifan}

Let  $G$ be a graph, $e=rs_1\in E(G)$ and $\varphi\in \CC^k(G-e)$ for some integer $k\ge 0$.
A \emph{multifan} centered at $r$ with respect to $e$ and $\varphi$
is a sequence $F_\varphi(r,s_1:s_p):=(r, rs_1, s_1, rs_2, s_2, \ldots, rs_p, s_p)$ with $p\geq 1$ consisting of  distinct vertices $r, s_1,s_2, \ldots , s_p$ and distinct edges $rs_1, rs_2,\ldots, rs_p$ satisfying the following condition:
\begin{enumerate}[(F1)]
	\item For every edge $rs_i$ with $i\in [2,p]$,  there exists $j\in [1,i-1]$ such that 
	$\varphi(rs_i)\in \pbar(s_j)$. 
\end{enumerate}
We will simply denote a multifan  $F_\varphi(r,s_1: s_{p})$ by $F$ if 
$\varphi$ and the vertices and edges in $F_\varphi(r,s_1: s_{p})$  are clear. 
Let $F_\varphi(r,s_1: s_{p})$ be a multifan. 
By its definition, for any $p^*\in [1,p]$,  $F_\varphi(r,s_1: s_{p^*})$
is a multifan. 
The following result regarding a multifan can be found in \cite[Theorem~2.1]{StiebSTF-Book}.

\begin{LEM}
	\label{thm:vizing-fan1}
	Let $G$ be a Class 2 graph and $F_\varphi(r,s_1:s_p)$  be a multifan with respect to a critical edge $e=rs_1$ and a coloring $\varphi\in \CC^\Delta(G-e)$. Then  the following statements  hold. 
	 \begin{enumerate}[(a)]
	 	\item $V(F)$ is $\varphi$-elementary. \label{thm:vizing-fan1a}
	 	\item Let $\alpha\in \pbar(r)$. Then for every $i\in [1,p]$  and $\beta\in \pbar(s_i)$,  $r$ 
	 	and $s_i$ are $(\alpha,\beta)$-linked with respect to $\varphi$. \label{thm:vizing-fan1b}
	 \end{enumerate}
\end{LEM}

Let $F_\varphi(r,s_1:s_p)$  be a multifan.  We call $s_{\ell_1},s_{\ell_2}, \ldots, s_{\ell_k}$, a subsequence of $s_1:s_p$, an  \emph{$\alpha$}-sequence with respect to $\varphi$ and $F$ if the following holds:
$$
\varphi(rs_{\ell_1})= \alpha\in \pbar(s_1),  \quad \varphi(rs_{\ell_i})\in \pbar(s_{\ell_{i-1}}), \quad  i\in [2,k].
$$
A vertex in an $\alpha$-sequence is called an \emph{$\alpha$-inducing vertex} with respect to $\varphi$ and $F$, and a missing color at an $\alpha$-inducing vertex is called an \emph{$\alpha$-inducing color}. For convenience, $\alpha$ itself is also an $\alpha$-inducing color. We say $\beta$ is {\it induced by} $\alpha$ if $\beta$ is $\alpha$-inducing. By Lemma~\ref{thm:vizing-fan1} (a) and the definition of multifan, each color in $\pbar(V(F))$ is induced by a unique color in $\pbar(s_1)$. Also if $\alpha_1,\alpha_2$ are two distinct colors in $\pbar(s_1)$, then an $\alpha_1$-sequence is disjoint with an $\alpha_2$-sequence. For two distinct $\alpha$-inducing colors $\beta$ and $\delta$, we write {$\mathit \delta \prec \beta$} if there exists an $\alpha$-sequence $s_{\ell_1},s_{\ell_2}, \ldots, s_{\ell_k}$ such that $\delta\in\pbar(s_{\ell_i})$, $\beta\in\pbar(s_{\ell_j})$ and $i<j$. For convenience, $\alpha\prec\beta$ for any $\alpha$-inducing color $\beta\not=\alpha$.
As a consequence of Lemma~\ref{thm:vizing-fan1} (a), we have the following properties for a multifan. 
A proof of the result can be found in~\cite[Lemma 3.2]{HZ}. 
\begin{LEM}
	\label{thm:vizing-fan2}
	Let $G$ be a Class 2 graph and $F_\varphi(r,s_1:s_p)$  be a multifan with respect to a critical edge $e=rs_1$ and a coloring $\varphi\in \CC^\Delta(G-e)$. For two colors $\delta\in \pbar(s_i)$ and $\lambda\in \pbar(s_j)$ with  $i,j\in [1,p]$ and $i\ne j$, the following statements  hold.
	\begin{enumerate}[(a)]
		\item If $\delta$ and $\lambda$ are induced by different colors, then $s_i$ and $s_j$ are $(\delta, \lambda)$-linked with respect to $\varphi$. 
		\label{thm:vizing-fan2-a}
		\item If $\delta$ and $\lambda$ are induced by the same color, $\delta\prec\lambda$, and $s_i$ and $s_j$ are $(\delta, \lambda)$-unlinked with respect to $\varphi$, 
		then $r\in P_{s_j}(\lambda, \delta, \varphi)$.  	\label{thm:vizing-fan2-b}
	\end{enumerate}
	
\end{LEM}

\subsection{Kierstead path}

Let $G$ be a graph, $e=v_0v_1\in E(G)$, and  $\varphi\in \CC^k(G-e)$ for some integer $k\ge 0$.
A \emph{Kierstead path}  with respect to $e$ and $\varphi$
is a sequence $K=(v_0, v_0v_1, v_1, v_1v_2, v_2, \ldots, v_{p-1}, v_{p-1}v_p,  v_p)$ with $p\geq 1$ consisting of  distinct vertices $v_0,v_1, \ldots , v_p$ and distinct edges $v_0v_1, v_1v_2,\ldots, v_{p-1}v_p$ satisfying the following condition:
\begin{enumerate}[(K1)]
	\item For every edge $v_{i-1}v_i$ with $i\in [2,p]$,  there exists $j\in [1,i-1]$ such that 
	$\varphi(v_{i-1}v_i)\in \pbar(v_j)$. 
\end{enumerate}

Clearly a Kierstead path with at most 3 vertices is a multifan. We consider Kierstead paths with $4$ vertices. The result below was proved in Theorem 3.3 from~\cite{StiebSTF-Book}. 

\begin{LEM}[]\label{Lemma:kierstead path1}
	Let $G$ be a Class 2 graph,
	 $e=v_0v_1\in E(G)$ be a critical edge, and $\varphi\in \CC^\Delta(G-e)$. If $K=(v_0, v_0v_1, v_1, v_1v_2,  v_2, v_2v_3, v_3)$ is a Kierstead path with respect to $e$
		and $\varphi$, then the following statements hold.
	\begin{enumerate}[(a)]
 		\item If $\min\{d_G(v_2), d_G(v_3)\}<\Delta$, then $V(K)$ is $\varphi$-elementary.
 		\item $|\pbar(v_3)\cap (\pbar(v_0)\cup \pbar(v_1))|\le 1$. 
 	\end{enumerate}

\end{LEM}

	\section{Proof of Theorems~\ref{thm:vertex spliting} and~\ref{thm:Delta-critical}}
We will prove Theorems~\ref{thm:vertex spliting} and~\ref{thm:Delta-critical}  based on the following lemmas, whose proof will be 
presented in the last section. 	

General properties on Kierstead paths with 5 vertices was proved by the first author 
of this paper~\cite{K5}. Here we stress only one of the cases. 
\begin{LEM}\label{lem:5vexKpathsettingup}
	Let $G$ be a Class 2 graph, $ab\in E(G)$ be a critical edge,  $\varphi\in \CC^\Delta(G-ab)$, and  
	$K=(a, ab,b,bu,u, us, s, st, t)$ be a Kierstead path with respect to $ab$ and $\varphi$. 
	If $|\pbar(t)\cap(\pbar(a)\cup \pbar(b))|\ge 3$, then the following hold:
	\begin{enumerate}[(a)]
		\item 
		There exists $\varphi^*\in \CC^\Delta(G-ab)$
		satisfies the following properties:
		\begin{enumerate}[(i)]
			\item  $\varphi^*(bu)\in \pbar^*(a)\cap \pbar^*(t)$, 
			\item  $\varphi^*(us)\in \pbar^*(b)\cap \pbar^*(t)$, and 
			\item $\varphi^*(st)\in \pbar^*(a)$. 
		\end{enumerate}
		\item  $d_G(b)=d_G(u)=\Delta$. 
		
	\end{enumerate}
	Figure~\ref{f3} shows a Kierstead path with the properties described in (a). 
\end{LEM}
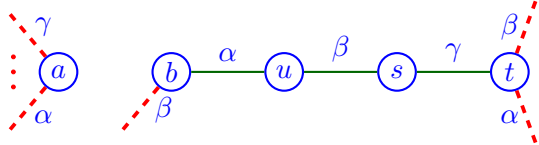
\begin{figure}[!htb]
	\begin{center}
		\begin{tikzpicture}[scale=1,rotate=90]
		
		{\tikzstyle{every node}=[draw ,circle,fill=white, minimum size=0.5cm,
			inner sep=0pt]
			\draw[blue,thick](0,-2) node (a)  {$a$};
			\draw[blue,thick](0,-3.5) node (b)  {$b$};
			
			\draw [blue,thick](0, -5) node (u)  {$u$};
			\draw [blue,thick](0, -6.5) node (s)  {$s$};
			\draw [blue,thick](0, -8) node (t)  {$t$};
		}
		\path[draw,thick,black!60!green]

		(b) edge node[name=la,pos=0.7, above] {\color{blue} $\alpha$\quad\quad} (u)
		(u) edge node[name=la,pos=0.7, above] {\color{blue}$\beta$\quad\quad} (s)
		(s) edge node[name=la,pos=0.7,above] {\color{blue} $\gamma$\quad\quad} (t);
		

		\draw[dashed, red, line width=0.5mm] (b)--++(140:1cm); 
		\draw[dashed, red, line width=0.5mm] (t)--++(200:1cm); 
		\draw[dashed, red, line width=0.5mm] (t)--++(340:1cm);
		\draw[dashed, red, line width=0.5mm] (a)--++(40:1cm); 
		
		\draw[dashed, red, line width=0.5mm] (a)--++(140:1cm);

		\draw[blue] (-0.5, -3.4 ) node {$\beta$};
		\draw[blue] (0.6, -1.8) node {$\gamma$};
		\draw[blue] (-0.6, -1.8) node {$\alpha$};
		\draw[blue] (0.6, -8.) node {$\beta$};
		\draw[blue] (-0.6, -8.) node {$\alpha$};
		
		{\tikzstyle{every node}=[draw, red ,circle,fill=red, minimum size=0.05cm,
			inner sep=0pt]
			\draw(-0.2,-1.4) node (f1)  {};
			\draw(0,-1.4) node (f1)  {};
			\draw(0.2, -1.4) node (f1)  {};
			
		} 
		
		\end{tikzpicture}

	\end{center}
	\caption{Colors on a Kierstead path of 5 vertices}
	\label{f3}
\end{figure}

 \begin{LEM}\label{lem:5vexKpathsettingup2}
	Let $G$ be a Class 2 graph, $ab\in E(G)$ be a critical edge,  $\varphi\in \CC^\Delta(G-ab)$, and  
	$K=(a, ab,b,bu,u, us, s, st, t)$  and $K^*=(a,ab,b,bu,u,ux,x)$ be two Kierstead paths with respect to $ab$ and $\varphi$, where $x\not\in V(K)$. 
	If $|\pbar(t)\cap(\pbar(a)\cup \pbar(b))|\ge 4$ and $\pbar(x) \subseteq \pbar(a)\cup \pbar(b)$, then   $d_G(x)=\Delta$. 
\end{LEM}	

A {\it short-kite}  $H$ is a graph with 
$$V(H)=\{a,b,c,u,x,y\} \quad \text{and}\quad E(H)=\{ab,ac,bu,cu,ux,uy\}.$$ 
The lemma below reveals some properties of a short-kite with specified colors on its edges.

\begin{LEM}\label{lemma:class2-with-fullDpair2}
	Let  $G$ be  a Class 2 graph, 
	 $H\subseteq G$ 
	be a short-kite with $V(H)=\{a,b,c,u,x,y\}$, and let $\varphi\in \CC^\Delta(G-ab)$. 
	Suppose $$K=(a,ab,b,bu,u,ux,x) \quad \text{and} \quad K^*=(b,ab,a,ac,c,cu,u,uy)$$
	are two Kierstead path with respect to $ab$ and $\varphi$.  
	If $\pbar(x)\cup \pbar(y)\subseteq \pbar(a)\cup \pbar(b)$,  then $\max\{d_G(x),d_G(y)\}=\Delta $. 
\end{LEM}

	A {\it kite}  $H$ is a graph with 
	$$V(H)=\{a,b,c,u,s_1,s_2,t_1,t_2\} \quad \text{and}\quad E(H)=\{ab,ac,bu,cu,us_1,us_2, s_1t_1,s_2t_2\}.$$ 
	The lemma below reveals some properties of a kite with specified colors on its edges. 
	\begin{LEM}\label{lem:kite}
		Let $G$ be a Class 2 graph, $H\subseteq G$ 
		be a kite with $V(H)=\{a,b,c,u,s_1,s_2,t_1,t_2\}$, and let $\varphi\in \CC^\Delta(G-ab)$. 
		Suppose $$K=(a,ab,b,bu,u,us_1, s_1,s_1t_1,t_1) \quad \text{and} \quad K^*=(b,ab,a,ac,c,cu,u,us_2, s_2,s_2t_2,t_2)$$
		are two Kierstead paths with respect to $ab$ and $\varphi$. 
		If $\varphi(s_1t_1)=\varphi(s_2t_2)$, 
		then $|\pbar(t_1)\cap \pbar(t_2) \cap ( \pbar(a)\cup \pbar(b))|\le 4$.  
	\end{LEM}

	Let $G$ be a $\Delta$-critical graph, $ab\in E(G)$, and $\varphi \in \CC^\Delta(G-ab)$. 
	A {\it fork}  $H$  with respect to $\varphi$ is a graph with 
	$$V(H)=\{a,b,u,s_1,s_2,t_1,t_2\} \quad \text{and}\quad E(H)=\{ab,bu,us_1,us_2, s_1t_1,s_2t_2\}$$ 
	such that $\varphi(bu)\in \pbar(a)$, $\varphi(us_1), \varphi(us_2) \in \pbar(a)\cup \pbar(b)$,  and $\varphi(s_1t_1)\in (\pbar(a)\cup \pbar(b))\cap \pbar(t_2) $
	and $\varphi(s_2t_2)\in (\pbar(a)\cup \pbar(b))\cap \pbar(t_1)$. 
	Fork  was defined in~\cite{av2} and  
	it was shown in~\cite[Proposition B]{av2} that a fork can not exist in a $\Delta$-critical graph if 
	the degree sum of $a$, $t_1$ and $t_1$ is small.

	\begin{LEM}\label{lem:fork}
	Let $G$ be a $\Delta$-critical graph, $ab\in E(G)$, and $\{u,s_1,s_2, t_1,t_2\}\subseteq V(G)$. 
	If $\Delta\ge d_G(a)+d_G(t_1)+d_G(t_2)+1$, 
	then for any $\varphi\in \CC^\Delta(G-ab)$, $G$ does not contain a fork on $\{a,b,u,s_1,s_2,t_1,t_2\}$ with respect to $\varphi$.   
\end{LEM}

We need the following two additional results to prove Theorem~\ref{thm:Delta-critical}. 
Since all vertices not missing  a given color $\alpha$
are saturated by the matching that consists of all edges colored by $\alpha$ in $G$, we have the 
following result.  

\begin{LEM}[Parity Lemma]
	Let $G$ be an $n$-vertex graph and $\varphi\in \CC^\Delta(G)$. 
	Then for any color $\alpha\in [1,\Delta]$, 
	$|\{v\in V(G): \alpha\in \pbar(v)\}| \equiv n \pmod{2}$. 
\end{LEM}

\begin{LEM}\label{lemma:class2-with-fullDpair}
	If  $G$ is an $n$-vertex Class 2 graph with a full-deficiency pair $(a,b)$ such that $ab$ is a critical edge of $G$, 
	then $G$ satisfies the following properties. 
	\begin{enumerate}[$(i)$]
		
		\item For every $x\in (N_G(a)\cup N_G(b))\setminus\{a,b\}$, $d_G(x)=\Delta$;
		
		\item For every $x\in V(G)\setminus\{a,b\}$, if $\dist_G(x, \{a,b\})=2$, then 
		$d_G(x)\ge \Delta-1$. 
		Furthermore, 
		if $d_G(a)<\Delta$ and $d_G(b)<\Delta$, then  $d_G(x)=\Delta$;
		\item  For  
		every   $x\in V(G)\setminus\{a,b\}$,  if $d_G(x)\ge n-|N_G(b)\cup N_G(a)|$, 
		then $d_G(x)\ge \Delta-1$. 
		Furthermore, 
		if $d_G(a)<\Delta$ and $d_G(b)<\Delta$, then  $d_G(x)=\Delta$;
		\item If there exists $x\in V(G)\setminus \{a,b\}$ such that $d_G(x)<\Delta$, then there exists 
		$y\in V(G)\setminus \{a,b,x\}$ such that $d_G(y)<\Delta$.
		%
	\end{enumerate} 
\end{LEM}
 \pf We let
 $\varphi\in \CC^\Delta(G-ab)$ and 
 $$
 F=(b, ba,a) 
 $$
 be the multifan  with respect to $ab$ and $\varphi$. 
 By Lemma~\ref{thm:vizing-fan1}~\eqref{thm:vizing-fan1a}, 
 \begin{equation}\label{pbarFa}
 |\pbar(F)|= 2\Delta+2-(d_{G}(a)+d_{G}(b))= 2\Delta+2-(\Delta+2)=\Delta.
 \end{equation}
 
 By Lemma~\ref{thm:vizing-fan1}, for every $\varphi'\in \CC^\Delta(G-ab)$,  $\{a,b\}$ is $\varphi'$-elementary and for every 
 $i\in \pbar'(a)$ and $j\in \pbar'(b)$, $a$ and $b$ are $(i,j)$-linked with respect to $\varphi'$. We will use this fact very often.

Since all the $\Delta$ colors appear in $\pbar(F)$,  each of $N_G(a)\cup \{b\}$ and $ N_G(b)\cup \{a\}$
is the vertex set of a  multifan with respect to $ab$ and $\varphi$. 
By Lemma~\ref{thm:vizing-fan1}~\eqref{thm:vizing-fan1a} and~\eqref{pbarFa}, 
we know that for every $x\in (N_G(a)\cup N_G(b))\setminus\{a,b\}$, $d_G(x)=\Delta$. 
This proves (i).

For (ii), let $x\in V(G)\setminus\{a,b\}$ such that $\dist_G(x, \{a,b\})=2$. 
We assume that $\dist_G(x, b)=2$ and let 
$u\in ( N_G(b))\setminus\{a\})\cap N_G(x)$.
Then by~\eqref{pbarFa}, $K=(a,ab, b, bu, u, ux, x)$ is a Kierstead path with respect to $ab$
and $\varphi$. By~\eqref{pbarFa} and Lemma~\ref{Lemma:kierstead path1} (b),  it follows that 
$d_G(x)\ge \Delta-1$. If $d_G(a)<\Delta$ and $d_G(b)<\Delta$, by~\eqref{pbarFa} and Lemma~\ref{Lemma:kierstead path1} (a), we get $d_G(x)=\Delta$.

For (iii), let $x\in V(G)\setminus\{a,b\}$ such that $d_G(x)\ge n-|N_G(b)\cup N_G(a)|$. By (i), we may assume that 
$x\not\in (N_G(a)\cup N_G(b))\setminus\{a,b\}$. Thus   $d_G(x)\ge n-|N_G(b)\cup N_G(a)|$ implies that 
there exists $u\in ( (N_G(a)\cup N_G(b)))\cap N_G(x)$. 
Therefore, $\dist_G(x,\{a,b\})=2$. Now Statement (ii) yields the conclusion. 
Statement (iv) is a consequence of~\eqref{pbarFa} and the Parity Lemma. 
%
\qed 

\begin{COR}\label{cor:no2D-1}
	Let  $G$ be  an $n$-vertex Class 2 graph with a full-deficiency pair $(a,b)$ such that $ab$ is a critical edge of $G$.  
	If $\Delta\ge \frac{3(n-1)}{4}$, then there exists at most one vertex $x\in V(G)\setminus \{a,b\}$ such that $d_G(x)=\Delta-1$.
\end{COR}	

\pf Assume to the contrary that there exist distinct $x,y\in V(G)\setminus \{a,b\}$ 
such that $d_G(x)=d_G(y)=\Delta-1$. By Lemma~\ref{lemma:class2-with-fullDpair} (i), $x,y\not\in (N_G(a)\cup N_G(b))\setminus\{a,b\}$. 
By Lemma~\ref{lemma:class2-with-fullDpair} (iii), we may assume  that $d_G(b)=\Delta$. 
Thus $d_G(a)=2$ as $d_G(a)+d_G(b)=\Delta+2$. 
Let $c$ be the other neighbor of $a$ in $G$. Since $(a,c)$ is a full-deficiency pair of $G$ as well,  we may assume $x,y\not\in N_G(c)$.

Since $d_{G}(b)=d_{G}(c)=\Delta$ and  $d_{G}(x)=d_{G}(y)=\Delta-1$, we get $|N_{G}(b)\cap N_{G}(c)|\ge \frac{n}{2}-1$ and $|N_{G}(x)\cap N_{G}(y)|\ge \frac{n}{2}-2$. 
Since $b,c,x,y\not\in N_{G}(b)\cap N_{G}(c)$ and $b,c,x,y\not\in N_{G}(x)\cap N_{G}(y)$, we get 
$|N_{G}(b)\cap N_{G}(c)\cap N_{G}(x)\cap N_G(y)|\ge 1$.
Let 
$u\in N_{G}(b)\cap N_{G}(c)\cap N_{G}(x)\cap N_{G}(y)$,  $H$ be 
the short-kite with $V(H)=\{a,b,c,u,x,y\}$, and let $\varphi\in \CC^\Delta(G-ab)$. 
As $\{a,b\}$
is $\varphi$-elementary$, |\pbar(a)\cup \pbar(b)|=2\Delta+2-(d_{G}(a)+d_{G}(b))= \Delta$ and so $\pbar(a)\cup \pbar(b)=[1,\Delta]$. Thus   
$K=(a,ab,b,bu,u,ux,x)$ and $ K^*=(b,ab,a,ac,c,cu,u,uy)$
are two Kierstead paths with respect to $ab$ and $\varphi$, and $\pbar(x)\cup \pbar(y)\subseteq \pbar(a)\cup \pbar(b)$.   
However,  $d_{G}(x)=d_{G}(y)=\Delta-1$,  contradicting  Lemma~\ref{lemma:class2-with-fullDpair2}.
\qed  

\proof[\bf Proof of Theorem~\ref{thm:vertex spliting}]
Since $G$ is overfull, $G$ is  Class 2. 
We show that every edge of $G$ is critical. 
Suppose to the contrary that there exists $xy\in E(G)$ such that $xy$ 
is not a critical edge of $G$.  Let 
$$
G^*=G-xy.  
$$
Then $\chi'(G^*)=\Delta+1$.

Since $ab$ is a critical edge of $G$,  $ab \ne xy$. 
Also, since $ab$ is a critical edge of $G$, and any $\Delta$-coloring of $G-ab$
gives a $\Delta$-coloring of $G^*-ab$, $ab$ is also a critical edge of $G^*$. 
Since $d_{G^*}(x)=d_{G^*}(y)=\Delta-1$, we reach a contradiction to Corollary~\ref{cor:no2D-1}. 
\qed

\proof[\bf Proof of Theorem~\ref{thm:Delta-critical}]  Let $(a,b)$ be a full-deficiency pair of $G$. 
It suffices to only show that for every $v\in V(G)\setminus\{a,b\}$, 
$d_G(v)=\Delta$.   To see this, let $G$ be a $\Delta$-critical graph with a
full-deficiency pair $(a,b)$ and  for every $v\in V(G)\setminus\{a,b\}$, 
$d_G(v)=\Delta$. Let $\varphi\in \CC^\Delta(G-ab)$. 
Since $\pbar(a)\cap \pbar(b)=\emptyset$ and $d_G(a)+d_G(b)=\Delta+2$, 
$\varphi(a)\cap \varphi(b)=\emptyset$. 
Thus, identifying $a$ and $b$ in $G$ gives a $\Delta$-coloring of a $\Delta$-regular multigraph $G^*$. 
This implies that $|V(G^*)|=n-1$ is even. 
So $n$ is odd. Consequently, $G$ is overfull.

Thus, for the sake of contradiction, we assume  
that there exists $x\in V(G)\setminus\{a,b\}$
such that $d_G(x)<\Delta$. By Lemma~\ref{lemma:class2-with-fullDpair} (iv), 
there exits $y\in V(G)\setminus\{a,b,x\}$ such that $d_G(y)<\Delta$. 
Furthermore, 
by Lemma~\ref{lemma:class2-with-fullDpair} (iii) and Corollary~\ref{cor:no2D-1},  there exists at most one vertex $x\in V(G)\setminus\{a,b\}$
such that $d_G(x)=\Delta-1$, and for all other vertex $y\in V(G)\setminus\{a,b,x\}$, if $d_G(y)<\Delta$,
then $d_G(y)<n-|N_G(b)\cup N_G(a)|$.  
This gives a vertex $t\in  V(G)\setminus\{a,b\}$ such that $d_G(t)<n-|N_G(b)\cup N_G(a)|$.  
Let $\varphi\in \CC^\Delta(G-ab)$ and 
$$
F=(b, ba,a) 
$$
be the multifan  with respect to $ab$ and $\varphi$. 
By Lemma~\ref{thm:vizing-fan1}~\eqref{thm:vizing-fan1a}, 
\begin{equation}\label{pbarF4}
|\pbar(F)|= 2\Delta+2-(d_{G}(a)+d_{G}(b))=\Delta.
\end{equation}


Assume, without loss of generality, that $d_G(b)\ge d_G(a)$. Then $d_G(b)\ge \frac{3(n-1)}{8}+1$  
as $d_G(a)+d_G(b)=\Delta+2\ge \frac{3(n-1)}{4}+2$.  
By Lemma~\ref{lemma:class2-with-fullDpair} (i) and (ii), we assume that $\dist_G(t, \{a,b\})\ge 3$. 
Since $\Delta\ge \frac{3(n-1)}{4}$,  for any $s\in N_G(t)$ with $d_G(s)\ge \Delta-1$ (such $s$ exists as $t$ is adjacent to at least two $\Delta$-neighbors by VAL),  we conclude that 
there exists $u\in N_G(b)\cap N_G(s)$.  Now by~\eqref{pbarF4}, $K=(a,ab,b,bu,u,us, s, st)$
is a Kierstaed path with respect to $ab$ and $\varphi$. 
This implies that $d_G(b)=d_G(u)=\Delta$ by  Lemma~\ref{lem:5vexKpathsettingup} (b).
Thus $d_G(a)=2$. We let $c$ be the other $\Delta$-neighbor of $a$. 

As $d_G(a)=2$ and $ab\in E(G)$, $|N_G(b)\cup N_G(a)| \ge \Delta+1> \frac{3n}{4}$. 
Since $G$ is $\Delta$-critical, VAL implies that for 
every $s\in N_G(t)$, $d_G(s)\ge \Delta+2-d_G(t)\ge \Delta+2 +|N_G(b)\cup N_G(a)|-n\ge n-|N_G(b)\cup N_G(a)|$. 
Thus, by Lemma~\ref{lemma:class2-with-fullDpair} (iii) and  Corollary~\ref{cor:no2D-1}, there exists at most 
one vertex $s\in N_G(t)$ such that $d_G(s)<\Delta$. In this case, $d_G(s)=\Delta-1$.

Next, we claim that 
\begin{equation}\label{claimx}
\text{for any $x\in V(G)\setminus \{a,b\}$,  $d_G(x)<\Delta \Rightarrow d_G(x)<n-|N_G(b)\cup N_G(a)|\le n-\Delta-1$.}
\end{equation}

Assume to the contrary that $d_G(x)\ge n-|N_G(b)\cup N_G(a)|$. By Lemma~\ref{lemma:class2-with-fullDpair} (iii), we have $d_G(x)=\Delta-1$. 
Again, as  $\Delta\ge \frac{3(n-1)}{4}$ and every vertex in $(N_G(a)\cup N_G(b)\cup N_G(t)) \setminus\{a,b\}$ has degree at least $\Delta-1$, for any $s\in N_G(t)$,  we conclude that 
there exists $u\in N_G(b)\cap N_G(s)\cap N_G(x)$.
Now by~\eqref{pbarF4}, $K=(a,ab,b,bu,u,us, s, st)$ and $K^*=(a,ab,b,bu,u,ux,x)$ 
are two Kierstaed paths with respect to $ab$ and $\varphi$. 
Clearly, $\pbar(t)\subseteq \pbar(a)\cup \pbar(b)$ and $\pbar(x)\subseteq \pbar(a)\cup \pbar(b)$,
and $|\pbar(t)|\ge \Delta-(n-\Delta-1)=2\Delta-n+1\ge \frac{n-1}{2}$.
Since $n\ge |V(K)\cup \{x\}|=6$,  $\Delta \ge  \frac{3(n-1)}{4}$ implies that $\Delta \ge 4$. 
As $\{b,s,t,x\}\cap N_G(b)=\emptyset$, we see that $n\ge 4+4=8 $. 
Hence, $|\pbar(t)|\ge \lceil \frac{n-1}{2} \rceil\ge 4$,  achieving a contradiction to Lemma~\ref{lem:5vexKpathsettingup2}.

By Lemma~\ref{lemma:class2-with-fullDpair} (iii) and (iv) and the conclusion in~\eqref{claimx},  
we let  $t_1,t_2\in V(G)\setminus \{a,b\}$ such that 
both of them have degree less than $n-\Delta$.  
Let $s_1\in N_G(t_1)$ and $s_2\in N_G(t_2)$ be any two distinct vertices. 
Since $G$ is $\Delta$-critical, VAL implies that for 
every $s_i\in N_G(t_i)$, $d_G(s_i)\ge 2\Delta-n>n-\Delta$. 
Thus, By Lemma~\ref{lemma:class2-with-fullDpair} (iii) and (iv) and the conclusion in~\eqref{claimx}, $d_G(s_1)=d_G(s_2)=\Delta$.   
Thus, since $b,c,s_1,s_2\not\in N_G(b)\cap N_G(c)$ and $b,c,s_1,s_2\not\in N_G(s_1)\cap N_G(s_2)$, 
\begin{eqnarray*}
	|N_G(b)\cap N_G(c)\cap N_G(s_1)\cap N_G(s_2)|&\ge& |N_G(s_1)\cap N_G(s_2)|-(n-|N_G(b)\cap N_G(c)|-4), \\
	&\ge & |N_G(s_1)\cap N_G(s_2)|+|N_G(b)\cap N_G(c)|+4-n\\
	&\ge & 2\Delta-n+2\Delta-n+4-n\ge 1,
\end{eqnarray*}
as $\Delta\ge \frac{3(n-1)}{4}$. Let $u\in N_G(b)\cap N_G(c)\cap N_G(s_1)\cap N_G(s_2)$. Then  $H$ with $V(H)=\{a,b,c,u,s_1,s_2,t_1,t_2\}$
is kite. By~\eqref{pbarF4}, both $$K=(a,ab,b,bu,u,us_1, s_1,s_1t_1,t_1) \quad \text{and} \quad K^*=(b,ab,a,ac,c,cu,u,us_2, s_2,s_2t_2,t_2)$$
are  Kierstead paths with respect to $ab$ and $\varphi$.
%
%
Let 
$$
A=\{t\in V(G)\setminus\{a,b\}: d_G(t)<n-\Delta \}. 
$$
We consider two cases below. 

\medskip 

{\bf \noindent Case 1: there exist two distinct $t_1,t_2\in A$ 
	such that $\varphi(t_1)\cap \varphi(t_2)\ne \emptyset$.}

\medskip 

In this case, we choose $s_1\in N_G(t_1)$ and $s_2\in N_G(t_2)$ such that $\varphi(s_1t_1)=\varphi(s_2t_2)$.
Let  $\Gamma= \pbar(t_1)\cap \pbar(t_2)$.  Since $\pbar(a)\cup \pbar(b)=[1,\Delta]$, $\Gamma\subseteq \pbar(a)\cup \pbar(b)$. 
By~\eqref{pbarF4} and the assumption of this case, 
$|\Gamma|\ge \Delta-2(n-\Delta-2)=3\Delta-2n+4\ge \frac{n-1}{4}+2$. 
Since $\Delta\ge \frac{3(n-1)}{4}\ge \frac{3}{4}|(V(H)|-1)$, $\Delta\ge 6$. 
Since also $s_1,s_2,t_1,t_2\not\in N_G(b)$, 
we have $n\ge |V(H)|+3\ge 11$. 
Thus,  $|\Gamma|\ge  \lceil \frac{n-1}{4}\rceil+2\ge 5$, 
contradicting  Lemma~\ref{lem:kite}.

\medskip 

{\bf \noindent Case 2: for each two distinct $t_1,t_2\in A$,  it holds that $\varphi(t_1)\cap \varphi(t_2)=\emptyset$.}

\medskip 
By~\eqref{pbarF4} and the assumption of this case, we see that $H^*$ with $V(H^*)=\{a,b,u,s_1,s_2,t_1,t_2\}$
is fork. However, by~\eqref{claimx},  $d_G(a)+d_G(t_1)+d_G(t_2)\le 2+2(n-\Delta-2)=2n-2\Delta-2 <\Delta$, 
as $\Delta\ge \frac{3(n-1)}{4}$, contradicting Lemma~\ref{lem:fork}. 
The proof is now completed.

\qed

\section{Proof of Lemmas~\ref{lem:5vexKpathsettingup} to~\ref{lem:kite}}


\begin{LEM1}
	Let $G$ be a Class 2 graph, $ab\in E(G)$ be a critical edge,  $\varphi\in \CC^\Delta(G-ab)$, and  
	 $K=(a, ab,b,bu,u, us, s, st, t)$ be a Kierstead path with respect to $ab$ and $\varphi$. 
	  If $|\pbar(t)\cap(\pbar(a)\cup \pbar(b))|\ge 3$, then the following hold:
	\begin{enumerate}[(a)]
		\item 
	There exists $\varphi^*\in \CC^\Delta(G-ab)$
	satisfies the following properties:
	\begin{enumerate}[(i)]
		\item  $\varphi^*(bu)\in \pbar^*(a)\cap \pbar^*(t)$, 
		\item  $\varphi^*(us)\in \pbar^*(b)\cap \pbar^*(t)$, and 
		\item $\varphi^*(st)\in \pbar^*(a)$. 
	\end{enumerate}
	\item  $d_G(b)=d_G(u)=\Delta$. 

		\end{enumerate}	
\end{LEM1}
\pf By Lemma~\ref{thm:vizing-fan1}, for every $\varphi'\in \CC^\Delta(G-ab)$,  $\{a,b\}$ is $\varphi'$-elementary and for every 
$i\in \pbar'(a)$ and $j\in \pbar'(b)$, $a$ and $b$ are $(i,j)$-linked with respect to $\varphi'$. 

Let $\Gamma=\pbar(t)\cap(\pbar(a)\cup \pbar(b))$, and 
 $\alpha,\beta \in \Gamma$. 
If $\alpha, \beta \in \pbar(a)$, then we let $\lambda\in \pbar(b)$, and do a $(\beta,\lambda)$-swap 
at $b$. If $\alpha, \beta \in \pbar(b)$, then we let $\lambda\in \pbar(a)$, and do a $(\beta,\lambda)$-swap 
at $a$.  Therefore, 
we may assume that $$\alpha\in \pbar(a) \quad \text{and} \quad \beta\in \pbar(b).$$

If $\varphi(bu)=\delta\ne \alpha$, then we do an $(\alpha,\delta)$-swap at $t$, 
and rename the color $\delta$ as $\alpha$ and vice versa. 
Thus we may assume $$\varphi(bu)=\alpha.$$

Assume first that $\varphi(us)\in \pbar(b)$. We do a $(\beta,\varphi(us))$-swap at $t$
and still call the resulting coloring by $\varphi$, we see that $\varphi(us)\in \pbar(b)\cap \pbar(t)$. 
By permuting the name of the colors, we let  $\varphi(us)=\beta$.
Let $\varphi(st)=\gamma$.  Since $\alpha,\beta\in \pbar(t)$, $\gamma\ne \alpha, \beta$. 
If $\gamma\in \pbar(a)$, we are done. So we assume  $\gamma\in \pbar(b)\cup \pbar(u)$.  
We color $ab$ by $\alpha$ and uncolor $bu$. Denote this resulting coloring by $\varphi'$.
Then $K'=(b,bu, u, us,s,st,t)$ is a Kierstead path with respect to $bu$ and $\varphi'$. 
However, $\alpha,\beta\in \pbar'(t)\cap (\pbar'(b)\cup \pbar'(u))$, showing a contradiction to 
Lemma~\ref{Lemma:kierstead path1} (b).

Thus we let $\varphi(us)=\delta\in \pbar(a)$. Then $\delta\ne \beta$
by Lemma~\ref{thm:vizing-fan1}~\eqref{thm:vizing-fan1a}.
Let $\varphi(st)=\gamma$. Clearly, $\gamma\ne \alpha,\beta, \delta$. We have either $\gamma\in \pbar(a)$
or $\gamma\in \pbar(b)\cup \pbar(u)$. 
We consider three cases below. 

{\bf \noindent Case 1. $\gamma\in \pbar(b)$.}

If $u\in P_a(\beta,\delta)=P_b(\beta,\delta)$,  we do a $(\beta,\delta)$-swap at $t$.
Since $a$ and $b$ are $(\delta,\gamma)$-linked and $u\in P_t(\delta,\gamma)$,  
we do a $(\delta,\gamma)$-swap at $a$. This gives a desired coloring $\varphi^*$. 

If $u\not\in P_a(\beta,\delta)=P_b(\beta,\delta)$, we first do a $(\beta,\delta)$-swap at $a$ and then a $(\beta,\gamma)$-swap at $a$. Again this gives a desired coloring $\varphi^*$.   

\smallskip 

{\bf \noindent Case 2. $\gamma\in \pbar(u)$.}

If $\delta\in \Gamma$, since $b$ and $u$
are $(\beta,\gamma)$-linked by Lemma~\ref{thm:vizing-fan1}~\eqref{thm:vizing-fan1b},  we do $(\beta,\gamma)$-swap at $t$. 
Now $u\in P_t(\delta,\beta)$, we  do a $(\beta,\delta)$-swap at $a$. 
This gives a desired coloring $\varphi^*$. 
Thus we assume $\delta\not\in \Gamma$. Since $b$ and $u$
are $(\beta,\gamma)$-linked by Lemma~\ref{thm:vizing-fan1}~\eqref{thm:vizing-fan1b}, 
and $a$ and $u$ are $(\delta,\gamma)$-linked by Lemma~\ref{thm:vizing-fan2}~\eqref{thm:vizing-fan2-a}, 
we do $(\beta,\gamma)-(\gamma,\delta)$-swaps at $t$. 
Finally, since $u\in P_t(\beta,\delta)$, we do a $(\beta,\delta)$-swap at $a$.
This gives a desired coloring $\varphi^*$.

\medskip 

{\bf \noindent Case 3. $\gamma\in \pbar(a)$.}

If $\delta\in \Gamma$, we do a $(\beta,\gamma)$-swap at $t$ and then a $(\beta,\delta)$-swap at $a$ to get a desired coloring $\varphi^*$. 
Thus we assume $\delta\not\in \Gamma$. Let $\tau\in \Gamma\setminus \{\alpha,\beta\}$. 
If $\tau\in \pbar(u)$, since $a$ and $u$ are $(\delta,\gamma)$-linked by Lemma~\ref{thm:vizing-fan2}~\eqref{thm:vizing-fan2-a}, we do a $(\tau,\delta)$-swap at $t$. This gives back to the previous case that 
$\delta\in \Gamma$. Next we assume $\tau\in \pbar(b)$. 
It is clear that $u\in P_a(\tau,\delta)=P_b(\tau,\delta)$, as otherwise, a $(\tau,\delta)$-swap at $a$ gives a desired coloring. 
Thus we do a $(\tau,\delta)$-swap at $t$,  giving back to the previous case that 
$\delta\in \Gamma$.

Now we assume  $\tau\in \pbar(a)$. 
If $u\not\in P_a(\beta,\delta)$, we do a $(\beta,\delta)$-swap at $a$. Since $a$ and $b$ are $(\alpha,\delta)$-linked and $u\in P_a(\alpha,\delta)$, we do an $(\alpha,\delta)$ swap at $t$.
Now since $u\in P_t(\gamma,\delta)$, we do a  
$(\gamma,\delta)$-swap at $a$, and do  $(\beta,\gamma)-(\gamma,\alpha)$-swaps at $t$.
Since $a$ and $b$ are $(\tau,\gamma)$-linked, we do a 
$(\tau,\gamma)$-swap at $t$, and then a $(\beta,\gamma)$
-swap at $a$. Now since $u\in P_t(\beta, \delta)$, 
we do a $(\beta, \delta)$-swap at $a$. This gives a desired coloring. 
Thus, we assume $u\in P_a(\beta,\delta)$. We do a $(\beta,\delta)$-swap at $t$, and then a $(\tau,\beta)$-swap at $t$.  
Next we do a $(\beta,\gamma)$-swap at $a$ and then a $(\gamma,\delta)$-swap at $a$. 
This gives a desired coloring for (b). 


For statement (b), let $\varphi^*\in \CC^\Delta(G-ab)$ satisfying (i)--(iii). 
Let $\alpha,\gamma\in \pbar^*(a)$, $\beta\in \pbar^*(b)$  with $\alpha,\beta \in \pbar(t)$ such that 
$$
\varphi^*(bu)=\alpha, \quad \varphi^*(us)=\beta,\quad \text{and} \quad  \varphi^*(st)=\gamma. 
$$
Let  $\tau\in \pbar^*(t)\setminus\{\alpha,\beta\}$. 
Suppose to the contrary first that $d_G(b)\le \Delta-1$. Let $\lambda\in \pbar^*(b)\setminus\{\beta\}$. 
We do $(\tau,\lambda)-(\lambda,\gamma)$-swaps at $t$. 
Now we color $ab$ by $\alpha$ and uncolor $bu$ to get a coloring $\varphi'$. 
Then $K'=(b,bu, u, us, s, st, t)$ is a Kierstead path with respect to $bu$ and $\varphi'$. However, $\alpha,\beta \in \pbar'(t)\cap (\pbar'(b)\cup \pbar'(u))$, contradicting Lemma~\ref{Lemma:kierstead path1} (b). 

Assume then that $ d_G(b)=\Delta$ and $d_G(u)\le \Delta-1$. Let $\lambda\in \pbar^*(u)$.
Since $(a,ab,b,bu,u)$ is a multifan, $\lambda\not\in \{\alpha,\beta,\gamma\}$.
Since $u$ and $b$ are $(\beta,\lambda)$-linked and $u$ and $a$ are 
$(\gamma,\lambda)$-linked by Lemma~\ref{thm:vizing-fan2}~\eqref{thm:vizing-fan2-b},  
we do $(\beta,\lambda)-(\lambda,\gamma)$-swap(s) at $t$. 
Now we color $ab$ by $\alpha$ and uncolor $bu$ to get a coloring $\varphi'$. 
Then $K'=(b,bu, u, us, s, st, t)$ is a Kierstead path with respect to $bu$ and $\varphi'$. However, $\alpha \in \pbar'(t)\cap \pbar'(u)$, contradicting Lemma~\ref{Lemma:kierstead path1} (a), since $d_G(u)<\Delta$.  
\qed

 \begin{LEM2}
 	Let $G$ be a Class 2 graph, $ab\in E(G)$ be a critical edge,  $\varphi\in \CC^\Delta(G-ab)$, and  
 	$K=(a, ab,b,bu,u, us, s, st, t)$  and $K^*=(a,ab,b,bu,u,ux,x)$ be two Kierstead paths with respect to $ab$ and $\varphi$, where $x\not\in V(K)$. 
 	If $|\pbar(t)\cap(\pbar(a)\cup \pbar(b))|\ge 4$ and $\pbar(x) \subseteq \pbar(a)\cup \pbar(b)$, then   $d_G(x)=\Delta$. 
\end{LEM2}	

\pf 
 Assume to the contrary that $d_G(x)\le \Delta-1$. 
 Since $\pbar(x) \subseteq \pbar(a)\cup \pbar(b)$,  
 Lemma~\ref{Lemma:kierstead path1} (b) gives that $d_G(x)= \Delta-1$. 
By Lemma~\ref{lem:5vexKpathsettingup}, $d_G(b)=d_G(u)=\Delta$ and we assume 
that $\pbar(b)=\beta$, 
$\varphi(bu)=\alpha$, $\varphi(us)=\beta$, $\varphi(st)=\gamma$, $\alpha,\gamma \in \pbar(a)$, 
and $\alpha,\beta\in \pbar(t)$. 
 In the following, when we swap colors, we always make sure that the 
  colors on the edges $bu$ and $us$ are unchanged. The color on the edge $st$ 
 might be changed, but the new color will still be a color from $\pbar(a)$. 
This is guaranteed by using the elementary fact that for every coloring $\varphi'\in \CC^\Delta(G-ab)$,
$a$ and $b$ are $(i,j)$-linked for every $i\in \pbar'(a)$ and every $j\in \pbar'(b)$.  
We use this fact every often without even mentioning it.   

Let $\varphi(ux)=\delta$ and $\pbar(x)=\tau$.   We first claim that if $\delta\ne \gamma$, then we may assume 
$\delta\in \pbar(t)$. 
Clearly, $\delta\ne \alpha,\beta$, and 
since $K^*$ is a Kierstead path and $\pbar(b)=\beta$,
we have $\varphi(ux)=\delta\in \pbar(a)$. 
Let $\Gamma=\pbar(t)\cap(\pbar(a)\cup \pbar(b))$, and let  
$\{\alpha,\beta,\eta,\lambda\}\subseteq \Gamma$. 
Suppose that $\delta\not\in \pbar(t)$. 
We do $(\beta,\gamma)-(\gamma,\eta)$-swaps at $b$. Denote the new coloring by $\varphi'$. 

If $ux\not\in P_b(\eta,\delta)$,  based on $\varphi'$, 
we do an $(\eta,\delta)$ -swap at $b$ and 
and then do an $(\alpha,\delta)$-swap 
at $t$. If $\tau=\gamma$, 
 we do $(\delta,\gamma)-(\gamma,\eta)$-swaps at $b$ and then do an $(\alpha,\eta)$-swap at $t$. 
Finally we do $(\eta,\lambda)-(\lambda,\gamma)-(\gamma,\beta)$-swaps at $b$. 
Thus we assume $\tau \ne \gamma$. Clearly, $\tau \ne \alpha$. 
If $\tau=\beta$, we simply do a $(\beta,\delta)$-swap at $b$, then $(\beta,\gamma)-(\gamma,\eta)$-swaps at $b$,  an $(\alpha,\eta)$-swap at $t$, and finally $(\eta,\lambda)-(\lambda,\gamma)-(\gamma,\beta)$-swaps at $b$. 
Thus, $\tau\ne \alpha,\beta$. We do $(\delta,\tau)-(\tau,\eta)$-swaps at $b$, an $(\alpha,\eta)$-swap at $t$, 
and finally do  $(\eta,\lambda)-(\lambda,\gamma)-(\gamma,\beta)$-swaps at $b$.

Thus, we assume that  $ux\not\in P_b(\eta,\delta)$. Based on $\varphi'$, we do an $(\eta,\delta)$-swap at $t$
and then $(\eta,\lambda)-(\lambda,\gamma)-(\gamma,\beta)$-swaps at $b$. 

After the operations above, we have $\varphi(bu)=\alpha$, $\varphi(us)=\beta$, $\varphi(st)=\gamma$,  $\varphi(ux)=\delta$ and $\pbar(x)=\tau$,  and $\alpha,\beta,\delta,\lambda \in \Gamma$.

\smallskip

{\noindent \bf Case 1: $\pbar(x)=\gamma$}.

\smallskip 

Recall that $\alpha,\beta,\delta\in \Gamma$. 
We color $ab$ by $\alpha$, recolor $bu$ by $\beta$, 
and uncolor $us$. 
Note that  $u$ and $t$
are $(\alpha,\gamma)$-linked, as otherwise an  $(\alpha,\gamma)$-swap 
at $u$ and a $(\beta,\gamma)$-swap at $s$ gives a coloring $\varphi'$
such that $\gamma\in \pbar'(u)\cap \pbar'(s)$. 
Thus we do an $(\alpha,\gamma)$-swap at both $a$ and $x$, recolor $ux$ by $\alpha$,  and then a $(\beta,\delta)$-swap 
at both $x$ and $a$. It is clear that $ux\in P_t(\alpha,\gamma)$ and $P_t(\alpha,\gamma)$ 
meets $u$ before $x$. We now do the following operations:
\[
\begin{bmatrix}
P_{[t,u]}(\alpha,\gamma)& ux & bu  &  ab  \\
\alpha/ \gamma& \alpha \rightarrow \beta &  \beta \rightarrow \gamma&
\gamma \rightarrow \beta \end{bmatrix}.
\]
Based on the coloring above, we do  $(\alpha,\beta)-(\beta,\delta)$-swaps at both $x$ and $a$,
and then an $(\alpha,\delta)$-swap at $x$. Denote the new coloring by $\varphi'$. 
Now we do an $(\alpha,\beta)$-swap at $t$
and color $us$ by $\alpha$, giving a $\Delta$-coloring of $G$. 

\smallskip 

{\noindent \bf Case 2: $\pbar(x)\ne \gamma$ and $\varphi(ux)\ne \varphi(st)=\gamma$}.

\begin{figure}[!htb]
	\begin{center}
		\begin{tikzpicture}[scale=1,rotate=90]
		
		{\tikzstyle{every node}=[draw ,circle,fill=white, minimum size=0.5cm,
			inner sep=0pt]
			\draw[blue,thick](0,-2) node (a)  {$a$};
			\draw[blue,thick](-1,-3) node (b)  {$b$};
			\draw [blue,thick](0, -5) node (u)  {$u$};
			\draw [blue,thick](-1, -7) node (x)  {$s$};
			\draw [blue,thick](-1, -9) node (t)  {$t$};
			\draw [blue,thick](1, -7) node (y)  {$x$};
		}
		\path[draw,thick,black!60!green]
		
		(b) edge node[name=la,pos=0.5, above] {\color{blue} $\alpha$\quad\quad} (u)
		(u) edge node[name=la,pos=0.6, above] {\color{blue}$\beta$\quad} (x)
		(u) edge node[name=la,pos=0.4,above] {\color{blue}  \quad$\delta$} (y)
		(x) edge node[name=la,pos=0.4,above] {\color{blue}  \quad$\gamma$} (t);
		

		\draw[dashed, red, line width=0.5mm] (b)--++(140:1cm); 
		\draw[dashed, red, line width=0.5mm] (t)--++(200:1cm); 
		\draw[dashed, red, line width=0.5mm] (t)--++(340:1cm); 
		\draw[dashed, red, line width=0.5mm] (y)--++(340:1cm);
		\draw[dashed, red, line width=0.5mm] (a)--++(40:1cm); 
			\draw[dashed, red, line width=0.5mm] (a)--++(90:1cm); 
		
		\draw[dashed, red, line width=0.5mm] (a)--++(140:1cm);

		\draw[blue] (-1.5, -9.5) node {$\alpha$}; 
			\draw[blue] (-.5, -9.5) node {$\delta$}; 
		\draw[blue] (1.5, -7.5) node {$\beta$}; 
		\draw[blue] (-1.2, -2.5) node {$\beta$};
		\draw[blue] (0.6, -1.8) node {$\delta$};
		\draw[blue] (0.2, -1.3) node {$\gamma$};
		\draw[blue] (-0.6, -1.8) node {$\alpha$};

%
%
		\end{tikzpicture}
		-	  	\end{center}
	\caption{Colors on the edges of $K$ and $K^*$}
	\label{pic1}
\end{figure}
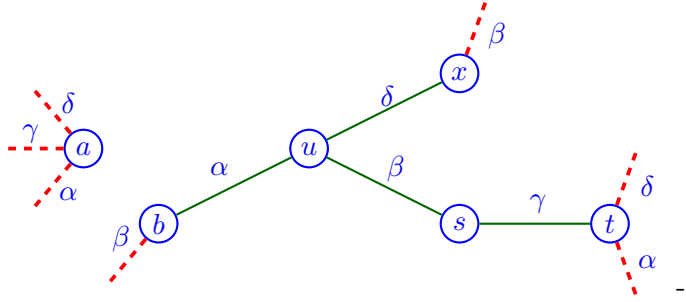

\smallskip 
Recall that $\alpha,\beta,\delta\in \Gamma$ and $|\Gamma|\ge 4$. 
Let  $\{\alpha,\beta,\delta,\lambda\}\subseteq \Gamma$. We show that there is a coloring $\varphi'$
such that $\{\alpha,\delta\}\subseteq \pbar'(t)$ and $\pbar'(x)=\beta$. 
Since we already have $\{\alpha,\delta\}\subseteq \pbar(t)$, we assume that $\pbar(x)=\tau\ne \beta$.  
Thus  $\tau\in \pbar(a)$. 
If $\tau=\alpha$, we simply do an $(\alpha,\beta)$-swap at $x$. 
Therefore, $\tau\ne \alpha$. It is possible that $\tau =\lambda$, but we deal with this 
together with the case that $\tau\ne \lambda$. 
We first do $(\beta,\gamma)-(\gamma,\lambda)-(\lambda,\tau)$-swaps at $b$, then we do an $(\alpha,\tau)$-swap at both $x$ and $t$. Now we do $(\tau,\gamma)-(\gamma,\beta)$-swaps at $b$,
and an $(\alpha,\beta)$-swap at both $x$ and $t$. 
We now derive a contradiction based on the coloring  of $E(K)\cup E(K^*)$, as shown in Figure~\ref{pic1}. 

We color $ab$ by $\alpha$, recolor $bu$ by $\beta$, 
and uncolor $us$. We then do an $(\alpha,\beta)$-swap at both $x$ and $t$, 
and an $(\alpha,\delta)$-swap at $x$. 
Now since $u$ and $s$ are $(\beta,\delta)$-linked, we do a $(\beta,\delta)$-swap at both $x$
and $a$.
Since $u$ and $t$ are $(\delta,\gamma)$-linked, we do a $(\delta,\gamma)$-swap at $a$. 
Finally, we do  $(\beta,\delta)-(\delta,\alpha)$-swaps at $x$. 
Now $P_u(\alpha,\beta)=uba$, so $u$ 
and $s$ are $(\alpha,\beta)$-unlined. We do an $(\alpha,\beta)$-swap at $u$
and color $us$ by $\beta$. 
This gives a $\Delta$-coloring of $G$, showing a contradiction. 

\smallskip 

{\noindent \bf Case 3: $\pbar(x)=\tau\ne \gamma$ and $\varphi(ux)=\varphi(st)=\gamma$}.

\smallskip 

We again let  $\Gamma=\pbar(t)\cap(\pbar(a)\cup \pbar(b))$, and let  
$\{\alpha,\beta,\lambda\}\subseteq \Gamma$.  We show that this case can be 
converted to Case 2. 
We first claim that $\pbar(x)=\tau\ne \beta$.
As otherwise, we first do an $(\alpha,\beta)$-swap at $x$,
and then a  $(\beta,\gamma)$-swap at $b$. 
Now, $P_b(\gamma,\alpha)=bux$, showing a contradiction to the fact that 
$a$ and $b$ are $(\alpha,\gamma)$-linked. 
Next, we claim that $\tau\ne \alpha$. 
As otherwise, we simply do a  $(\beta,\gamma)$-swap at $b$
and achieve a same contradiction as above. 
Thus, $\tau\ne \alpha,\beta$. 

We do a  $(\beta,\gamma)$-swap at $b$
and an $(\alpha,\gamma)$-swap at $t$. 
Now do a $(\tau,\gamma)$-swap at both $x$ and $t$, a $(\gamma,\lambda)$-swap at $b$,
an $(\lambda,\alpha)$-swap at $t$, and 
finally a 
$(\beta,\lambda)$-swap at $b$. Let the new coloring be $\varphi'$. 
We see that $\varphi'(st)=\lambda \ne \varphi'(ux)=\tau$. 
We verify that it still holds that $|\pbar'(t)\cap(\pbar'(a)\cup \pbar'(b))|\ge 4$. 
If $\tau \in \Gamma$, we now have $\alpha,\beta,\gamma,\tau \in \pbar'(t)\cap(\pbar'(a)\cup \pbar'(b))$.
If $\tau \not\in \Gamma$, then  $(\Gamma\setminus\{\lambda\} )\cup \{\tau\} \subseteq \pbar'(t)\cap(\pbar'(a)\cup \pbar'(b))$.
\qed 
 
\begin{LEM3}
	Let  $G$ be  a Class 2 graph, 
	$H\subseteq G$ 
	be a short-kite with $V(H)=\{a,b,c,u,x,y\}$, and let $\varphi\in \CC^\Delta(G-ab)$. 
	Suppose $$K=(a,ab,b,bu,u,ux,x) \quad \text{and} \quad K^*=(b,ab,a,ac,c,cu,u,uy)$$
	are two Kierstead path with respect to $ab$ and $\varphi$.  
	If $\pbar(x)\cup \pbar(y)\subseteq \pbar(a)\cup \pbar(b)$,  then $\max\{d_G(x),d_G(y)\}=\Delta $. 
\end{LEM3}
 
 \pf 
 Assume to the contrary that  $\max\{d_G(x),d_G(y)\}\le \Delta-1$. 
 Since both $K$ and $K^*$ are Kierstead paths and 
 $\pbar(x)\cup \pbar(y)\subseteq \pbar(a)\cup \pbar(b)$, Lemma~\ref{Lemma:kierstead path1} (a)
 and (b) implies that $d_G(b)=d_G(u)=\Delta$
 and $d_G(x)=d_G(y)=\Delta-1$. 
 
 
 Let $\pbar(b)=\{1\}$. Then $\varphi(ac)=1$. 
 We may assume  $\varphi(uy)=1$.  The reasoning is below. 
 Since $a$ and $b$ are $(1,\alpha)$-linked for every $\alpha\in \pbar(y)\subseteq \pbar(a)\cup \pbar(b)$, 
 we may assume  $\pbar(y)=1$. Then a $(1,\varphi(uy))$-swap at  $y$
 gives a coloring, call it still $\varphi$, such that $\varphi(uy)=1$. 
 We consider now two cases. 
 
 \smallskip 
 {\bf \noindent Case 1: $\pbar(x)=\pbar(y)$.} 
 \smallskip 
 
 Let 
 $ \varphi(ux)=\gamma$,  \text{and}  $\pbar(x)=\pbar(y)=\eta.$	
 As $\varphi(uy)=\pbar(b)=1$, $1\not\in \{\gamma, \eta\}$. 
 As both $K$ and $K^*$ are Kierstead paths and 
 $\pbar(x)\cup \pbar(y)\subseteq \pbar(a)\cup \pbar(b)$, $\gamma,\eta \in \pbar(a)$. 
 Denote by $P_u(1,\gamma)$  the $(1,\gamma)$-subchain starting at $u$ that does not include the edge $ux$. 
 \begin{CLA}\label{cla:claim1}
 	We may assume that $P_u(1,\gamma)$    ends at $x$, some vertex  $z\in V(G)\setminus\{a,b,c,u,x,y\}$,  or passing $c$ ends at $a$. 
 \end{CLA}
 
 \pf Note that $P_a(1,\gamma)=P_b(1,\gamma)$. If $u\not\in P_a(1,\gamma)$, then the $(1,\gamma)$-chain containing $u$ is a cycle or a path with endvertices contained in $V(G)\setminus\{a,b,c,u,x,y\}$.  Thus  
 $P_u(1,\gamma)$ ends at  $x$ or some  $z\in V(G)\setminus\{a,b,c,u,x,y\}$.  Hence we assume $u\in P_a(1,\gamma)$. 
 As a consequence,  $P_u(1,\gamma)$ ends at either $b$ or $a$. 
 If $P_x(1,\gamma)$ ends at $b$, we color $ab$ by 1, uncolor $ac$, and exchange the vertex labels $b$ and $c$. 
 This gives an edge $\Delta$-coloring of $G-ab$ such that  $P_u(1,\gamma)$ ends at  $a$.
 Thus, if $u\in P_a(1,\gamma)$, we may always assume that $P_u(1,\gamma)$ ends at  $a$.	  
 \qed 
 
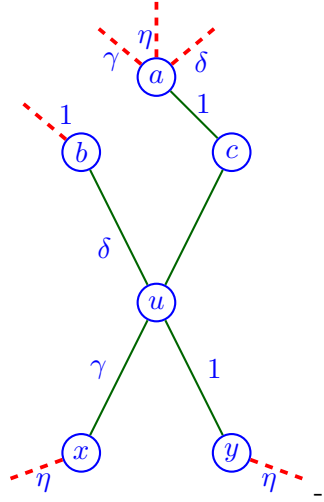
\begin{figure}[!htb]
 	\begin{center}
 		\begin{tikzpicture}[scale=1]
 		
 		{\tikzstyle{every node}=[draw ,circle,fill=white, minimum size=0.5cm,
 			inner sep=0pt]
 			\draw[blue,thick](0,-2) node (a)  {$a$};
 			\draw[blue,thick](-1,-3) node (b)  {$b$};
 			\draw[blue,thick](1,-3) node (c)  {$c$};
 			\draw [blue,thick](0, -5) node (u)  {$u$};
 			\draw [blue,thick](-1, -7) node (x)  {$x$};
 			\draw [blue,thick](1, -7) node (y)  {$y$};
 		}
 		\path[draw,thick,black!60!green]
 		(a) edge node[name=la,pos=0.7, above] {\color{blue} $1$} (c)
 		
 		(c) edge node[name=la,pos=0.5, below] {\color{blue}} (u)
 		(b) edge node[name=la,pos=0.5, below] {\color{blue} $\delta$\quad\quad} (u)
 		(u) edge node[name=la,pos=0.6, above] {\color{blue}$\gamma$\quad\quad} (x)
 		(u) edge node[name=la,pos=0.6,above] {\color{blue}  \quad$1$} (y);
 		

 		\draw[dashed, red, line width=0.5mm] (b)--++(140:1cm); 
 		\draw[dashed, red, line width=0.5mm] (x)--++(200:1cm); 
 		\draw[dashed, red, line width=0.5mm] (y)--++(340:1cm);
 		\draw[dashed, red, line width=0.5mm] (a)--++(40:1cm); 
 		\draw[dashed, red, line width=0.5mm] (a)--++(90:1cm); 
 		\draw[dashed, red, line width=0.5mm] (a)--++(140:1cm);

 		\draw[blue] (-1.5, -7.4) node {$\eta$}; 
 		\draw[blue] (1.5, -7.4) node {$\eta$}; 
 		\draw[blue] (-1.2, -2.5) node {$1$};
 		\draw[blue] (0.6, -1.8) node {$\delta$};
 		\draw[blue] (-0.6, -1.8) node {$\gamma$};
 		\draw[blue] (-0.15, -1.5) node {$\eta$};
 		
%
%
%
 		\end{tikzpicture}
 		-	  	\end{center}
 	\caption{Colors on the edges connecting $x$ and $y$ to $b$}
 	\label{f1}
 \end{figure}

 Let $\varphi(bu)=\delta$.  Again, $\delta\in \pbar(a)$. 
 Figure~\ref{f1} depicts the colors and missing colors  on these specified edges and vertices, respectively.  Clearly, $\delta\ne 1, \gamma$.
 Since $a$ and $b$ are $(1,\delta)$-linked with respect to $\varphi$, $\eta\ne \delta$. 
 Thus, $\gamma, \delta$ and $\eta$ are pairwise distinct.   
 
 \begin{CLA}\label{cla:claim2}
 	It holds that $ub\in P_y(\eta,\delta)$ and $P_y(\eta,\delta)$ meets $u$ before $b$. 
 \end{CLA}
 
 \pf Let $\varphi'$ be obtained from $\varphi$ by coloring $ab$ by $\delta$ and uncoloring $bu$. Note that $\pbar'(b)=1, \pbar'(u)=\delta$ and $\varphi'(uy)=1$.
 Thus $F^*=(u, ub,b, uy, y)$ is a multifan and so $u$ and $y$ are $(\eta, \delta)$-linked by Lemma~\ref{thm:vizing-fan1}~\eqref{thm:vizing-fan1b}. By uncoloring $ab$ and coloring 
 $bu$ by $\delta$, we get back the original coloring $\varphi$. Therefore, under the coloring $\varphi$, $u\in P_y(\eta, \delta)$ and $P_y(\eta,\delta)$ meets $u$ before $b$. 
 \qed 
 
 We apply the following operations based on $\varphi$:
 \[
 \begin{bmatrix}
 ux& P_{[u,y]}(\eta,\delta)  &  ub &  P_u(1,\gamma)  & ab\\
 \gamma\rightarrow \eta& \delta/\eta &  \delta \rightarrow 1&
 1/\gamma&  \delta\end{bmatrix}.
 \]
 By Claim~\ref{cla:claim1},  $P_u(1,\gamma)$ does not end at  $b$. 
 In any case, the above operations give  
 an edge $\Delta$-coloring of $G$. This contradicts  the earlier assumption that $\chi'(G)=\Delta+1$.
 
 \medskip 
 
 {\bf \noindent Case 2: $\pbar(x)\ne \pbar(y)$.}
 
 Let 
 $$ \varphi(bu)=\alpha, \quad  \varphi(ux)=\beta,  \quad \pbar(x)=\tau, \quad \text{and} \quad \pbar(y)=\gamma.$$	
 As $\varphi(uy)=\pbar(b)=1$, $1\not\in \{\alpha,\beta,\gamma\}$.
 Also, since $a$ and $b$ are $(1,\alpha)$-linked, $\gamma\ne \alpha$. 
 Since both $K$ and $K^*$ are Kierstead paths and 
 $\pbar(x)\cup \pbar(y)\subseteq \pbar(a)\cup \pbar(b)$, we have $\alpha,\beta,\tau,\gamma\in\pbar(a)$. 
 \begin{CLA}
 	We may assume $\pbar(x)=\tau=1$. 
 \end{CLA} 
 \pf  If $uy\not\in P_x(1,\tau)$, we simply do a $(1,\tau)$-swap at $x$. 
 Thus, we assume that $u\in P_x(1,\tau)$. We first do a $(1,\tau)$-swap at $b$, then an $(\alpha,\tau)$-swap at $x$. Then we do a $(\gamma,\tau)$-swap at $b$. Finally, a $(1,\gamma)$-swap at $b$ and a $(1,\alpha)$-swap at $x$
 give the desired coloring. 
 \qed 
 
\begin{figure}[!htb]
 	\begin{center}
 		\begin{tikzpicture}[scale=1]
 		
 		{\tikzstyle{every node}=[draw ,circle,fill=white, minimum size=0.5cm,
 			inner sep=0pt]
 			\draw[blue,thick](0,-2) node (a)  {$a$};
 			\draw[blue,thick](-1,-3) node (b)  {$b$};
 			\draw[blue,thick](1,-3) node (c)  {$c$};
 			\draw [blue,thick](0, -5) node (u)  {$u$};
 			\draw [blue,thick](-1, -7) node (x)  {$x$};
 			\draw [blue,thick](1, -7) node (y)  {$y$};
 		}
 		\path[draw,thick,black!60!green]
 		(a) edge node[name=la,pos=0.7, above] {\color{blue} $\beta$} (c)
 		
 		(c) edge node[name=la,pos=0.5, below] {\color{blue}} (u)
 		(b) edge node[name=la,pos=0.5, below] {\color{blue} $\alpha$\quad\quad} (u)
 		(u) edge node[name=la,pos=0.6, above] {\color{blue}$\beta$\quad\quad} (x)
 		(u) edge node[name=la,pos=0.6,above] {\color{blue}  \quad$1$} (y);
 		

 		\draw[dashed, red, line width=0.5mm] (b)--++(140:1cm); 
 		\draw[dashed, red, line width=0.5mm] (x)--++(200:1cm); 
 		\draw[dashed, red, line width=0.5mm] (y)--++(340:1cm);
 		\draw[dashed, red, line width=0.5mm] (a)--++(40:1cm); 
 			\draw[dashed, red, line width=0.5mm] (a)--++(90:1cm); 
 		\draw[dashed, red, line width=0.5mm] (a)--++(140:1cm);

 		\draw[blue] (-1.5, -7.4) node {$1$}; 
 		\draw[blue] (1.5, -7.4) node {$\gamma$}; 
 		\draw[blue] (-1.2, -2.5) node {$\beta$};
 		\draw[blue] (0.6, -1.8) node {$\gamma$};
 		\draw[blue] (-0.6, -1.8) node {$1$};
 			\draw[blue] (-0.15, -1.5) node {$\alpha$};

%

 		\draw [orange,thick](3.5, -5) node (t)  {$\Rightarrow$};

 		\begin{scope}[shift={(7,0)}]
 		{\tikzstyle{every node}=[draw ,circle,fill=white, minimum size=0.5cm,
 			inner sep=0pt]
 			\draw[blue,thick](0,-2) node (a)  {$a$};
 			\draw[blue,thick](-1,-3) node (b)  {$b$};
 			\draw[blue,thick](1,-3) node (c)  {$c$};
 			\draw [blue,thick](0, -5) node (u)  {$u$};
 			\draw [blue,thick](-1, -7) node (x)  {$x$};
 			\draw [blue,thick](1, -7) node (y)  {$y$};
 		}
 		\path[draw,thick,black!60!green]
 		(a) edge node[name=la,pos=0.7, above] {\color{blue} $\beta$} (c)
 		(a) edge node[name=la,pos=0.7, above] {\color{red} $\alpha$} (b)
 		
 		(c) edge node[name=la,pos=0.5, below] {\color{blue}} (u)
 		(b) edge node[name=la,pos=0.5, below] {\color{red} $\beta$\quad\quad} (u)
 		(u) edge node[name=la,pos=0.6,above] {\color{blue}  \quad$1$} (y);
 		

 		\draw[dashed, red, line width=0.5mm] (u)--++(200:1cm); 
 		\draw[dashed, red, line width=0.5mm] (x)--++(200:1cm); 
 		\draw[dashed, red, line width=0.5mm] (x)--++(340:1cm);
 		\draw[dashed, red, line width=0.5mm] (y)--++(340:1cm);
 		\draw[dashed, red, line width=0.5mm] (a)--++(40:1cm); 
 		\draw[dashed, red, line width=0.5mm] (a)--++(140:1cm);

 		\draw[blue] (-1.5, -7.4) node {$1$}; 
 		\draw[blue] (-0.5, -7.4) node {$\beta$}; 
 		\draw[blue] (1.5, -7.4) node {$\gamma$}; 
 		\draw[blue] (0.6, -1.8) node {$\gamma$};
 		\draw[blue] (-0.6, -1.8) node {$1$};
 		\draw[blue] (-0.5, -5.5) node {$\alpha$}; 
 			\draw[blue] (-0.5, -5.5) node {$\alpha$};

%
 		\end{scope}
 		\end{tikzpicture}
 		-	  	\end{center}
 	\caption{Colors on the edges connecting $x$ and $y$ to $b$}
 	\label{f2}
 \end{figure}
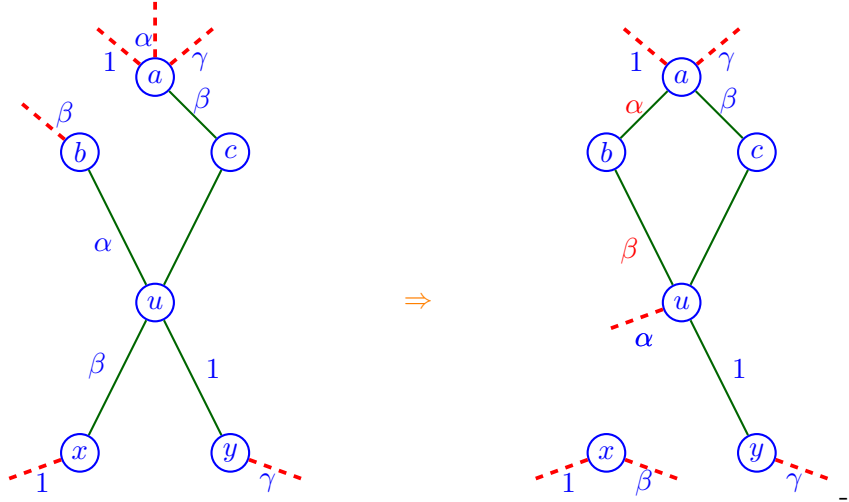

 Since $ux\in P_x(1,\beta)$, and $a$ and $b$ are $(1,\beta)$-linked, we do a $(1,\beta)$-swap at $b$. 
 Now we color $ab$ by $\alpha$, recolor $bu$ by $\beta$
 and uncolor $ux$, see Figure~\ref{f2} for a depiction. 
 
 Note that 
 $$
 F^*=(u, ux,x,uy,y), \quad K^*=(x,xu,u,ub,b,ba,a)
 $$
 are, respectively,  a multifan and 
 a Kierstead path. 
 By Lemma~\ref{thm:vizing-fan1}~\eqref{thm:vizing-fan1b}, $u$ and $y$ are $(\alpha,\gamma)$-linked, and $u$
 and $x$ are $(\alpha,\beta)$-linked and $(1,\alpha)$-linked. 
 Thus, we do an $(\alpha,\gamma)$-swap at $a$, an $(\alpha,\beta)$-swap at $a$, a $(1,\alpha)$-swap at $a$, 
 and then an $(\alpha,\gamma)$-swap at $a$. Now $P_u(\alpha,\beta)=uba$, contradicting 
 Lemma~\ref{thm:vizing-fan1}~\eqref{thm:vizing-fan1b} that $u$ and $x$ are $(\alpha,\beta)$-linked. 
 The proof  is now completed. 
 \qed

 \begin{LEM4}
 	Let $G$ be a Class 2 graph, $H\subseteq G$ 
 	be a kite with $V(H)=\{a,b,c,u,s_1,s_2,t_1,t_2\}$, and let $\varphi\in \CC^\Delta(G-ab)$. 
 	Suppose $$K=(a,ab,b,bu,u,us_1, s_1,s_1t_1,t_1) \quad \text{and} \quad K^*=(b,ab,a,ac,c,cu,u,us_2, s_2,s_2t_2,t_2)$$
 	are two Kierstead paths with respect to $ab$ and $\varphi$.
 	If $\varphi(s_1t_1)=\varphi(s_2t_2)$, 
 	then $|\pbar(t_1)\cap \pbar(t_2) \cap ( \pbar(a)\cup \pbar(b))|\le 4$.  
 \end{LEM4}
 
 \pf Let $\Gamma=\pbar(t_1)\cap \pbar(t_2) \cap ( \pbar(a)\cup \pbar(b))$. 
 Assume to the contrary that $|\Gamma|\ge 5$. By considering $K$ and applying Lemma~\ref{lem:5vexKpathsettingup}, we conclude that $d_G(b)=d_G(u)=\Delta$.
 We show that there exists $\varphi^*\in \CC^\Delta(G-ab)$ satisfying the following properties:
 \begin{enumerate}[(i)]
 	\item  $\varphi^*(bu), \varphi^*(cu), \varphi^*(us_2)\in \pbar^*(a)\cap \pbar^*(t_1)\cap \varphi^*(t_2)$, 
 	\item  $\varphi^*(us_1)\in \pbar^*(b)\cap \pbar^*(t_1)\cap \varphi^*(t_2)$, and 
 	\item $\varphi^*(s_1t_1)=\varphi^*(s_2t_2)\in \pbar^*(a)$. 
 \end{enumerate}

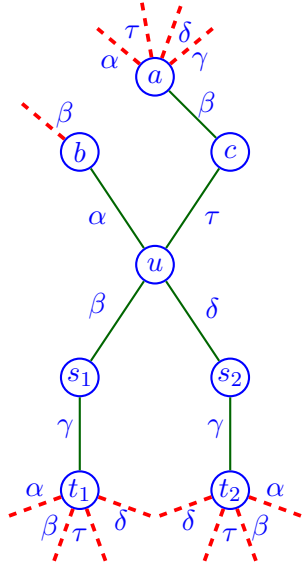
\begin{figure}[!htb]
	\begin{center}
		\begin{tikzpicture}[scale=1]
		
		{\tikzstyle{every node}=[draw ,circle,fill=white, minimum size=0.5cm,
			inner sep=0pt]
			\draw[blue,thick](0,-2) node (a)  {$a$};
			\draw[blue,thick](-1,-3) node (b)  {$b$};
			\draw[blue,thick](1,-3) node (c)  {$c$};
			\draw [blue,thick](0, -4.5) node (u)  {$u$};
			\draw [blue,thick](-1, -6) node (x)  {$s_1$};
			\draw [blue,thick](1, -6) node (y)  {$s_2$};
				\draw [blue,thick](-1, -7.5) node (t1)  {$t_1$};
			\draw [blue,thick](1, -7.5) node (t2)  {$t_2$};
		}
		\path[draw,thick,black!60!green]
		(a) edge node[name=la,pos=0.8, above] {\color{blue} $\beta$} (c)
		
		(c) edge node[name=la,pos=0.4, below] {\color{blue} \quad$\tau$} (u)
		(b) edge node[name=la,pos=0.4, below] {\color{blue} $\alpha$\quad\quad} (u)
		(u) edge node[name=la,pos=0.6, above] {\color{blue}$\beta$\quad\quad} (x)
			(x) edge node[name=la,pos=0.7, above] {\color{blue}$\gamma$\quad\quad} (t1)
		(y) edge node[name=la,pos=0.7,above] {\color{blue}  $\gamma$\quad\quad} (t2)
		(u) edge node[name=la,pos=0.6,above] {\color{blue}  \quad$\delta$} (y);
		

		\draw[dashed, red, line width=0.5mm] (b)--++(140:1cm); 
		\draw[dashed, red, line width=0.5mm] (t1)--++(200:1cm); 
		\draw[dashed, red, line width=0.5mm] (t1)--++(250:1cm); 
		\draw[dashed, red, line width=0.5mm] (t1)--++(290:1cm); 
		\draw[dashed, red, line width=0.5mm] (t1)--++(340:1cm); 
		\draw[dashed, red, line width=0.5mm] (t2)--++(200:1cm); 
		\draw[dashed, red, line width=0.5mm] (t2)--++(250:1cm); 
		\draw[dashed, red, line width=0.5mm] (t2)--++(290:1cm); 
		\draw[dashed, red, line width=0.5mm] (t2)--++(340:1cm); 
		
		\draw[dashed, red, line width=0.5mm] (a)--++(40:1cm); 
			\draw[dashed, red, line width=0.5mm] (a)--++(100:1cm); 
				\draw[dashed, red, line width=0.5mm] (a)--++(70:1cm); 
				\draw[dashed, red, line width=0.5mm] (a)--++(140:1cm); 

		\draw[blue] (-1.6, -9+1.5) node {$\alpha$}; 
		\draw[blue] (1.6, -9+1.5) node {$\alpha$}; 
		\draw[blue] (-1.4, -9.5+1.5) node {$\beta$}; 
		\draw[blue] (1.4, -9.5+1.5) node {$\beta$}; 
		\draw[blue] (-1, -9.6+1.5) node {$\tau$}; 
		\draw[blue] (1, -9.6+1.5) node {$\tau$}; 
		\draw[blue] (-0.45, -9.4+1.5) node {$\delta$}; 
		\draw[blue] (0.45, -9.4+1.5) node {$\delta$};

		\draw[blue] (-1.2, -2.5) node {$\beta$};
		\draw[blue] (0.6, -1.8) node {$\gamma$};
		\draw[blue] (-0.6, -1.8) node {$\alpha$};
			\draw[blue] (-0.3, -1.4) node {$\tau$};
				\draw[blue] (0.4, -1.4) node {$\delta$};
		
%
%
		
		\end{tikzpicture}
		  	\end{center}
	\caption{Colors on the edges of a kite}
	\label{pic2}
\end{figure}

See Figure~\ref{pic2} for a depiction of the colors described above. 

 Let $\alpha,\beta, \tau,\delta\in \Gamma$, and let $\varphi(s_1t_1)=\varphi(s_2t_2)=\gamma$. 
 We may assume that $\alpha\in \pbar(a)$ and $\beta\in \pbar(b)$. 
 Otherwise, since $d_G(b)=\Delta$, we have $\alpha,\beta \in \pbar(a)$. 
 Let $\lambda\in \pbar(b)$. As $a$ and $b$ are  $(\beta,\lambda)$-linked, we do a  
 $(\beta,\lambda)$-swap at $b$. Note that this operation may change some colors of the edges of $K$
 and $K^*$, but they are still Kierstead paths with respect to $ab$
 and the current coloring. 
 
 Since $d_G(b)=d_G(u)=\Delta$, and $\beta\in \pbar(b)\cap \pbar(t_1)$, 
 we know that $\gamma\in \pbar(a)$, as  $K_1$ is a
 Kierstead path. 
 Next, we may assume that $\varphi(bu)=\alpha$. 
 If not, let $\varphi(bu)=\alpha'$. Since $a$
 and $b$ are $(\alpha,\beta)$-linked, we do an $(\alpha,\beta)$-swap at $b$.
 Now $a$ and $b$ are   $(\alpha,\alpha')$-linked, we do an  $(\alpha,\alpha')$-swap at $b$. 
 Finally, we do an $(\alpha',\beta)$-swap at $b$. 
 All these swaps do not change the colors in $\Gamma$, so now we get the color on $bu$
 to be $\alpha$. 
 
 We may now assume that $\varphi(cu)=\tau$. 
 If not, let $\varphi(cu)=\tau'$. 
 Since $a$ and $b$ are $(\beta,\tau)$-linked, we do a $(\beta,\tau)$-swap at $b$. 
 Then do  $(\tau,\tau')-(\tau',\beta)$-swaps at $b$. 
 
 Finally, we show that we can modify $\varphi$
 to get $\varphi'$ such that $\varphi'(us_1)=\beta$
 and $\varphi'(us_2)=\delta$. 
 Assume firstly that $\varphi(us_1)=\beta'\ne \beta$. 
 If $\beta'\in \Gamma$, we do $(\beta,\gamma)-(\gamma,\beta')$-swaps at $b$. 
 Thus, we assume $\beta'\not\in \Gamma$. 
 Let $\lambda\in \Gamma\setminus \{\alpha,\beta,\tau,\delta\}$. 
 If $u\not\in P_a(\beta,\beta')=P_b(\beta,\beta')$, we simply do a $(\beta,\beta')$-swap at $b$. 
 Thus, we assume $u\in P_a(\beta,\beta')=P_b(\beta,\beta')$. 
 We do a $(\beta,\beta')$-swap at both $t_1$
 and $t_2$. Since $a$
 and $b$ are $(\beta,\lambda)$-linked, we 
 do a $(\beta,\lambda)$-swap at both $t_1$
 and $t_2$. 
 Now we do $(\beta,\gamma)-(\gamma,\beta')$-swaps at $b$. 
 By switching the role of $\beta$ and $\beta'$, 
 we have $\varphi(us_1)=\beta$. 
 Lastly, we show that $\varphi(us_2)=\delta$. 
 
 Note that $bu\in P_{t_1}(\alpha,\gamma)$. 
 Otherwise, let $\varphi'=\varphi/P_{t_1}(\alpha,\gamma)$.
 Then $P_b(\alpha,\beta)=bus_1t_1$, showing 
 a contradiction to the fact that $a$
 and $b$ are $(\alpha,\beta)$-linked with respect to $\varphi'$. 
 Thus, $bu\in P_{t_1}(\alpha,\gamma)$.
 Next, we claim that $P_{t_1}(\alpha,\gamma)$ meets $u$ before $b$. 
 As otherwise, we do the following operations to get a $\Delta$-coloring of $G$:
 \[
 \begin{bmatrix}
 s_1t_1& P_{[s_1,b]}(\alpha,\gamma)  &  us_1 &  bu  & ab\\
 \gamma\rightarrow \alpha& \alpha/\gamma &  \beta \rightarrow \alpha&
 \alpha\rightarrow \beta&  \gamma\end{bmatrix}.
 \]
 This gives a contradiction to the assumption that $G$ is $\Delta$-critical. 
 Thus, we  have that $P_{t_1}(\alpha,\gamma)$ meets $u$ before $b$.
 This implies that it is not the case that $P_{t_2}(\alpha,\gamma)$
 meets $u$ before $b$. In turn, this implies that $u\in P_a(\beta,\delta')=P_b(\beta,\delta')$.  
 As otherwise,  we get a $\Delta$-coloring of $G$ by doing a $(\beta,\delta)$-swap along the $(\beta,\delta)$-chain containing $u$,
 and then doing the same operation as above with $t_2$
 playing the role of $t_1$. 
 
 Since $u\in P_a(\beta,\delta')=P_b(\beta,\delta')$, 
 we do a $(\beta,\delta')$-swap at both $t_1$ and $t_2$. As $u\in P_a(\beta,\tau)=P_b(\beta,\tau)$, 
 we do a $(\beta,\tau)$-swap at both $t_1$
 and $t_2$. Since $us_1\in P_{t_1}(\beta,\gamma)$, 
 we do a $(\beta,\gamma)$-swap at $b$, then a  $(\gamma,\lambda)$-swap at $b$. 
 Since $a$ and $b$ are $(\tau,\lambda)$-linked, we do a $(\tau,\lambda)$-swap at both $t_1$
 and $t_2$. Now $(\lambda,\delta)-(\delta,\gamma)-(\gamma,\beta)$-swaps at $b$
 give a desired coloring. 
 
 Still, by the same arguments as above,  we have that  $P_{t_1}(\alpha,\gamma)$ meets $u$ before $b$,
 and   $u\in P_a(\beta,\delta)=P_b(\beta,\delta)$. 
 Let $P_u(\beta,\delta)$ be the $(\beta,\delta)$-chain starting at $u$ not including the edge $us_2$. 
 It is clear that $P_u(\beta,\delta)$  ends at either $a$ or $b$. 
 We may assume that $P_u(\beta,\delta)$ ends at $a$. 
 Otherwise, we color $ab$ by $\beta$, uncolor $ac$, and let $\tau$
 play the role of $\alpha$. Let $P_u(\alpha,\gamma)$
 be the $(\alpha,\gamma)$-chain starting at $u$ not including the edge $bu$,
 which ends at $t_1$ by our earlier argument. 
 We do the following operations to get a $\Delta$-coloring of $G$:
 \[
 \begin{bmatrix}
 P_u(\alpha,\gamma)& bu &  P_u(\beta,\delta) &  us_2t_2 & ab\\
 \alpha/\gamma& \alpha\rightarrow \beta &  \beta/\delta&
 \delta/\gamma&  \alpha\end{bmatrix}.
 \]
 This gives a contradiction to the assumption that $G$ is $\Delta$-critical.
The proof is now finished. 
 \qed

\bibliographystyle{plain}
\bibliography{SSL-BIB_08-19}
\end{document}